%% file: ks-symbDynII-main.tex
\documentclass[10pt]{article}
\usepackage{amssymb}
\usepackage{amsmath}
\usepackage{graphicx}
\usepackage{xspace}

\newtheorem{theorem}{Theorem}
\newtheorem{definition}{Definition}
\newtheorem{lemma}[theorem]{Lemma}

\newtheorem{remark}[theorem]{Remark}

\newcommand{\dom }{{\rm dom}\,}
\newcommand{\ltwoq}{l_{2,q}}
\newcommand{\cover}[1]{\stackrel{#1}{\Longrightarrow}}
\newcommand{\WqS}[2]{W_{#1,#2}}
\newcommand{\PM}{\mathcal P}
\newcommand{\gss}{GSS\xspace}
\newcommand{\VL}{{\rm \textbf{VC}}\xspace}
\newcommand{\NC}[1]{\textbf{N#1}}

\def\qed{{\hfill{\vrule height5pt width3pt depth0pt}\medskip}}
\def\comment#1{{}}

\begin{document}
\begin{center}
{\Large \bf   Symbolic dynamics for the Kuramoto-Sivashinsky PDE on the line: connecting orbits between periodic orbits}

 \vskip 0.5cm
{\large Daniel Wilczak} and {\large Piotr Zgliczy\'nski}\footnote{Work of D.W. and P.Z. was supported by National Science Center (NCN) of Poland under project No. UMO-2016/22/A/ST1/00077} % Maestro
 \vskip 0.2cm
  Jagiellonian University, Faculty of Mathematics and Computer Science, \\
 \L ojasiewicza 6, 30--348  Krak\'ow, Poland \\
 \texttt{e-mail:\,\{Daniel.Wilczak,Piotr.Zgliczynski\}@uj.edu.pl}
\vskip 0.5cm

\today
\end{center}

\vskip 0.5cm

\begin{abstract}
We prove the existence of infinite number of homoclinic and heteroclinic orbits to two periodic orbits for the Kuramoto-Sivashinsky PDE on the line with odd and periodic boundary conditions and for some fixed parameter value of the system. The proof is computer assisted and it is based on a new algorithm for rigorous integration of the variational equation for a class of dissipative PDEs on the torus.
\end{abstract}

\input intro.tex

\input top-tools.tex

\input applications.tex

\input varsysintro.tex

\input c1alg.tex

\appendix

%Appendix contains three parts.
% The first
%part consisting from Appendices \ref{sec:lognorm}, \ref{sec:gss}, \ref{sec:limitlogn}, \ref{sec:c1converproof} has a theoretical character is  an abstract exposition of our approach together with proofs of convergence theorems stated in Section~\ref{sec:c1-conver}. Appendix~\ref{sec:ConCondKs} contains a proof that for the sets used in our computations for KS equation all conditions required for our approach are satisfied.   The last part contains the explicit formulas for various objects appearing in our algorithm.

%\input lognorm.tex

%\input ks_convergence_cond.tex

\input appendix-num.tex

\input ref.tex
\end{document}

%% file: intro.tex
\section{Introduction}

The goal of the present work is twofold. First, we present a new algorithm for rigorous integration of the variational equation (i.e. producing $\mathcal C^1$ estimates) for a class of dissipative PDEs on the torus. For finite dimensional ODEs the $\mathcal C^1$ algorithm \cite{Z02} proved to be a very useful and efficient tool in analysis of various phenomena. The spectrum of applications is wide -- from proving existence and checking stability of periodic orbits \cite{KS,WB}, the existence of connecting orbits \cite{WZ03,W}, solving boundary value problems, proving existence of normally hyperbolic manifolds, validation of the existence of Arnold's diffusion far from perturbative case \cite{CW} to studying global behavior such as existence of hyperbolic attractors \cite{Tucker1,W10}. We believe our new $\mathcal C^1$ algorithm will provide a general tool for analysis of similar phenomena in a class of dissipative PDEs.

The second goal of this paper is to strengthen the results from \cite{WZ} about symbolic dynamics for the Kuramoto-Sivashinsky PDE on the line. In fact, such application was the main motivation for us to develop the $\mathcal C^1$ algorithm for PDEs. Using this new algorithm we proved the existence of countable infinity of homoclinic and heteroclinic orbits between two periodic orbits in the Kuramoto-Sivashinsky PDE \cite{KT,S} on the line -- see Theorem~\ref{thm:connectingOrbits} below.

We consider the one-dimensional Kuramoto-Sivashinsky PDE (in the sequel we will refer to it as the KS equation), which is given by
\begin{equation}\label{eq:KS}
u_t = -\nu u_{xxxx} - u_{xx} + (u^2)_x, \qquad \nu>0,
\end{equation}
where $x \in \mathbb{R}$, $u(t,x) \in \mathbb{R}$ and we impose
odd and periodic boundary conditions
\begin{equation}
u(t,x)=-u(t,-x), \qquad u(t,x) = u(t,x+ 2\pi). \label{eq:KSbc}
\end{equation}

The Kuramoto-Sivashinsky equation has been introduced by Kuramoto~\cite{KT} in space dimension one for the study of front propagation in the Belousov-Zhabotinsky reactions. An extension of this equation to space dimension 2 (or more)
has been given by G. Sivashinsky \cite{S} in studying the propagation of flame front in the case of mild combustion.

The following theorem  is the main result of \cite{WZ}.

\begin{theorem}\cite[Thm.~1]{WZ}\label{thm:symdynKS}
    The system (\ref{eq:KS})--(\ref{eq:KSbc}) with the parameter value $\nu=0.1212$ is chaotic in the following sense.
    There exists a compact invariant set $\mathcal A\subset L^2((-\pi,\pi))$ ($\mathcal{A}$ is compact in $H^k((-\pi,\pi))$ for any $k \in \mathbb{N}$) which consists of
    \begin{enumerate}
        \item bounded full trajectories visiting explicitly given and disjoint vicinities  of two selected periodic solutions $u^1$ and $u^2$, respectively, with any prescribed order $\{u^1,u^2\}^\mathbb Z$.% This gives rise to symbolic dynamics.
        \item countable infinity of periodic orbits with arbitrary large periods. In fact, each periodic sequence of symbols $\{u^1,u^2\}^\mathbb Z$ is realised by a periodic solutions of the system (\ref{eq:KS})--(\ref{eq:KSbc}).
    \end{enumerate}
\end{theorem}
\begin{figure}
    \centerline{
        \includegraphics[width=.5\textwidth]{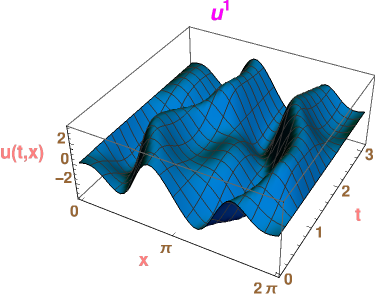}
        \includegraphics[width=.5\textwidth]{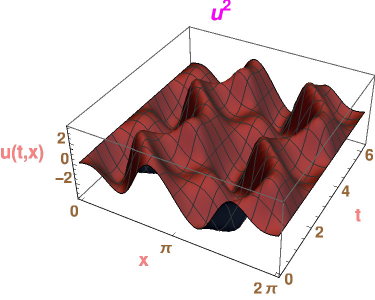}
    }
    \caption{Two approximate time-periodic orbits $u^1$ and $u^2$.}
    \label{fig:periodic-orbits-3d}
\end{figure}
The two special solutions $u^1$ and $u^2$ appearing in Theorem~\ref{thm:symdynKS} are time-periodic
--- see Fig.\ref{fig:periodic-orbits-3d}. In the present work  we extend this result and prove that the invariant set contains orbits asymptotic to $u^1$ and $u^2$ both in backward and forward time.

\begin{theorem}\label{thm:connectingOrbits}
    The invariant set $\mathcal A$ from Theorem~\ref{thm:symdynKS} contains countable infinity of geometrically different
    \begin{enumerate}
        \item homoclinic orbits to $u^1$,
        \item homoclinic orbits to $u^2$,
        \item heteroclinic orbits between $u^1$ and $u^2$,
        \item heteroclinic orbits between $u^2$ and $u^1$.
    \end{enumerate}
    The convergence in the definition of homoclinic and heteroclinic orbits happens in $L^2((-\pi,\pi))$ and in $H^k((-\pi,\pi))$ for any $k \in \mathbb{N}$.
\end{theorem}
Connecting orbits (homoclinic and heteroclinic) between periodic orbits $u^i$ and $u^j$ are understood as non-periodic solutions defined for all $t\in\mathbb R$, such that their $\alpha$-limit and $\omega$-limit sets in $L^2$ topology (turns out that also $H^k$ for any $k$) are equal $u^i$ and $u^j$, respectively. In fact, each non-constant sequence of symbols $(v_1,\ldots,v_n)\in \{u^1,u^2\}^n$ encodes a connecting orbit in the following way. There is a trajectory of the system (\ref{eq:KS})--(\ref{eq:KSbc}) with $\alpha$-limit set equal to $v_1$ and $\omega$-limit set equal to $v_n$. Moreover, the solution visits disjoint neighbourhoods of some points on $u^1$ and $u^2$ with the order prescribed by the sequence $(v_1,\ldots,v_n)$. Hence, we have countable infinity of geometrically different connecting orbits between $u^i$ and $u^j$, $i,j=1,2$.

  Our approach uses no special features of the Kuramoto-Sivashinsky PDE. We just need the Fourier basis to diagonalize the linear part of the PDE under consideration. Therefore, our method should be applicable to other systems of dissipative PDEs with periodic boundary conditions. While there is a wide literature on computer-assisted verification of symbolic dynamics for ODEs (see for example \cite{AZ,GZ,MM1,MM2,Mos,RNS,Tucker1,Tucker2,W,WSB,WZ03,Z4}), to the best of our knowledge, this is the first computer-assisted result of this type for PDEs.

The proofs of Theorems~\ref{thm:symdynKS} and~\ref{thm:connectingOrbits}  are a mixture of geometrical methods and rigorous numerics.
The proof of Theorem~\ref{thm:connectingOrbits} combines topological arguments used to construct symbolic dynamics (as in Theorem~\ref{thm:symdynKS}) with suitable cone conditions to obtain desired homoclinic and heteroclinic orbits to some periodic orbits.

In the case of the KS equation, we apply the method to certain Poincar\'e map. The orbits $u^1$ and $u^2$ from Theorem~\ref{thm:symdynKS} correspond then to a fixed point and a period two point. First, we construct symbolic dynamics along apparent heteroclinic connections between these points. Then using the cone conditions we show that these heteoroclinic indeed exist.  This is the same type of construction as it is used in the proof of the Smale-Birkhoff homoclinic theorem \cite[Thm.5.3.5]{GH} about the existence of horseshoe
dynamics for finite-dimensional diffeomorphism with a hyperbolic fixed point whose stable and unstable manifolds intersect transversally.

The nature of   our geometric method is such that its assumptions   can be expressed as a finite set of explicit inequalities. Therefore, we can use a computer to verify rigorously that these inequalities are satisfied for a given map $f$, provided we have an algorithm that computes rigorous bounds on $f$ on compact sets. In the case of the KS equation we have to compute rigorously bounds on some Poincar\'e maps and its partial derivatives.
For this purpose we propose an algorithm, which allows to compute rigorous bounds on the trajectories of PDEs with periodic boundary conditions
(this was described in \cite{WZ}) together with bounds on the derivatives with the respect to initial conditions.

Let us recall from \cite{WZ} the motivation for the choice of the KS equation for this study.
The following is known about the dynamics of KS on the line:
\begin{itemize}
\item The existence of a compact global attractor and the existence of a finite-dimensional inertial manifolds  for (\ref{eq:KS})--(\ref{eq:KSbc}) are well
established, see  \cite{CEES,FT,FNST,NST} and the
literature cited there. We would like to emphasise, that we are not
using these results in our work.
\item There exist multiple numerical studies of the
dynamics of the KS equation (see for example \cite{CCP,HN,JKT,JJK,SP}), where it
was shown, that the dynamics of the KS equation can be  highly
nontrivial for some values of parameter $\nu$, while being  well
represented by relatively small number of modes.
\item There are  several papers devoted to computer-assisted proofs of periodic
orbits for the KS equation by Zgliczy\'nski \cite{ZKSper,ZKS3}, Arioli and Koch
\cite{AK10} and by Figueras, Gameiro, de la Llave and Lessard \cite{FGLL,FL17,GL17}.
\end{itemize}

While the choice of the odd periodic boundary conditions was motivated by earlier
numerical studies of KS equation \cite{CCP,JKT},  the basic mathematical
reason is the following: the equation (\ref{eq:KS}) with periodic boundary
conditions has the translational symmetry. This implies, that for
a fixed value of $\nu$, all periodic orbits are members of one-parameter
families of periodic orbits. The restriction to the invariant subspace of odd functions breaks this symmetry, and gives a hope, that dynamically interesting objects are isolated and easier accessible for  computer-assisted proofs.

\begin{figure}
\centerline{\includegraphics[width=.8\textwidth]{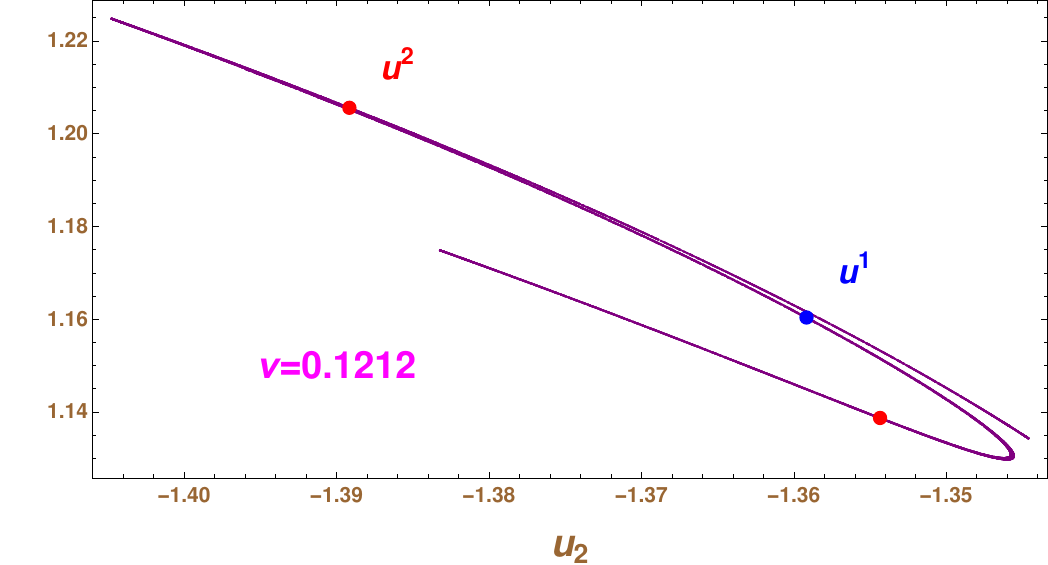}}
\caption{Numerically observed chaotic attractor for (\ref{eq:KS})--(\ref{eq:KSbc}) obtained by simulation of a finite-dimensional projection of the corresponding infinite-dimensional ODE for the Fourier coefficients in $u(t,x)=\sum_{k=1}^\infty u_k(t)\sin (kx)$. Projection onto $(u_2,u_3)$ plane of the intersection of the observed attractor with the Poincar\'e section $u_1=0, u_1'>0$ is shown along with an approximate location of the two periodic points $u^1, u^2$ appearing in Theorems~\ref{thm:symdynKS} and~\ref{thm:connectingOrbits}. The point $u^1$ is a fixed point for the Poincar\'e map and $u^2$ is of period two.}
\label{fig:attractor}
\end{figure}

Our choice of the parameter value $\nu=0.1212$ is motivated by a numerical  observation, that the Feigenbaum route to chaos through successive period doubling bifurcations \cite{F78}  happens for (\ref{eq:KS}) as $\nu$ decreases toward $\nu=0.1212$  (first observed for other values of $\nu$ in \cite{SP}). For this parameter value a~chaotic attractor is observed --- see Fig.~\ref{fig:attractor}.

Another result illustrating the effectiveness of our integration algorithm is about the attracting periodic orbit for $\nu=0.127$.
\begin{theorem}\label{thm:attrPerOrbit}
 Consider  system (\ref{eq:KS})--(\ref{eq:KSbc}) with the parameter value $\nu=0.127$ and  $\varphi(t,p)$ be the induced local semiflow.

 Then there exists a point $x_0$ and a compact forward invariant set $\mathcal{A}\subset H^4((-\pi,\pi))$, $x_0 \in \mathcal{A}$ % and $U$ - an open neighborhood of $x_0$ in %$L^2((-\pi,\pi))$ \textbf{TODO: $H^k((-\pi,\pi))$},
 such that
 \begin{itemize}
 %  \item for every $p \in U$, there exists $t>0$ such that $\varphi(t,p) \in \mathcal{A}$
   \item $x_0$ is time periodic
   \item   there exists $\lambda>0$ and $C>0$, such that for each $p_1,p_2 \in \mathcal{A}$ and $t>0$ holds
   \begin{equation*}
     \|\varphi(t,p_1) - \varphi(t,p_2)\|_\ast \leq Ce^{-\lambda t}\|p_1 - p_2\|_{\ast},
   \end{equation*}
   where $\| \cdot \|_\ast$ is some norm constructed using Fourier series (see section~\ref{sec:attr-per-orb} for details).
 \end{itemize}
\end{theorem}

The existence of this periodic orbit was first proved with computer assistance in \cite{ZKSper}
using the first version of the algorithm of integration of dissipative PDEs based on self-consistent bounds. In a sense it is one of the easiest
periodic orbits to prove, however its attracting character was not proved, yet.  In  papers \cite{ZKS3,AK10,FL17,GL17} the existence of multiple periodic
orbits for this system has been proved for more difficult parameter values, in \cite{AK10,GL17} the authors were able to show that some of these
orbits are unstable. It turns out that checking the instability of an orbit is an easier task than proving that the orbit is attracting.
In  proofs of instability in \cite{GL17} it was enough to prove the existence of one unstable direction, while the stability proof requires control in all directions in the phase space. This is what we do in the present paper using our $\mathcal C^1$-algorithm. With such an algorithm proving instability or stability of apparently hyperbolic orbits is the task of the same difficulty. Observe that in the proof of the instability of some periodic orbit in \cite{AK10} the whole spectrum has been estimated and a first $\mathcal C^1$ algorithm for KS-equation has been proposed.

The content of the paper can be described as follows. We have two parts: geometric and rigorous-numerics ones. The geometric part consists of Section~\ref{sec:geometricTools} containing a description of some geometric tools from dynamics, which are used in the proof of Theorem~\ref{thm:connectingOrbits}. In Section~\ref{sec:conn-orb-ks} computer-assisted proofs of Theorem~\ref{thm:connectingOrbits} and Theorem~\ref{thm:attrPerOrbit} are given. These proofs rely on rigorous estimates for some Poincar\'e maps induced by the KS equation. The rest of the paper, the second part, describes how we can obtain such bounds. In Section~\ref{sec:variationalEquations} we recall from \cite{WZ24} some results about convergence solutions of Galerkin projections of dissipative PDEs together with their derivatives with respect to initial conditions.
In Section~\ref{sec:c1algorithm} we present our $\mathcal{C}^1$-algorithm for integration of dissipative PDEs.

The main body of paper is followed by the appendix, which contains the explicit formulas for various objects appearing in our algorithm for the KS equation.

%% file: top-tools.tex
\section{Geometric tools from dynamics - covering relations and cone conditions}\label{sec:geometricTools}

The goal of this section is to describe the tools from dynamics, which we use to obtain the existence of hyperbolic periodic orbits and heteroclinic connections between them for the KS equation.

\subsection{Covering relations.}
%\label{subsec:infdim-covrel}
Topological methods proved to be very useful in the context of computer-assisted study of dynamical systems. One of the most efficient in the context of studying chaotic dynamics is the \emph{method of covering relations} introduced in \cite{Z0,Z4} for maps with one exit ("unstable") direction and later extended to include many exit directions in \cite{ZGi}  known also in the literature  as the method of correctly aligned windows \cite{E1,E2}. In \cite{WZ} this method has been adopted to the case of convex and compact subsets of the real normed spaces. Here we recall definition and the main topological result from \cite{WZ}.

\begin{definition}%\label{def:hset}
    Let $X$ be a real normed space. A $h$-set, $N=(|N|,c_N,u(N))$, is an object consisting of the following data
    \begin{itemize}
        \item $|N|\subset X$ -- a compact set called \emph{the support} of $N,$
        \item $u(N)$ -- a non-negative integer,
        \item $c_N$ -- a homeomorphism of $X$ such that
        \begin{equation*}
        c_N^{-1}(|N|) = \overline{B_{u(N)}}\oplus T_N =: N_c,
        \end{equation*}
        where $T_N$ is a convex set.
    \end{itemize}
\end{definition}
The above definition generalizes the concept of $h$-sets introduced in \cite{ZGi} for finite-dimensional spaces, where a $h$-set is defined as the product of closed unit balls $\overline{B}_u\times \overline{B}_s\subset \mathbb R^{u+s}$ in a coordinate system $c_N$.  The next definition extends the notion of covering relations from \cite{ZGi} to compact maps acting on infinite-dimensional real normed spaces.

\begin{definition}%\label{def:covrel}
    Let $N$ and $M$ be h-sets in $X$ and $Y$, respectively such that $u=u(N)=u(M)$. Let $f\colon |N|\to Y$ be a continuous map and set $f_c:=c_M^{-1} \circ f\circ c_N$. We say that $N$ $f$-covers $M$, denoted by $N\cover{f}M$,  if there is a linear map $L:\mathbb R^u\to \mathbb R^u$ and a compact homotopy $H:[0,1]\times N_c\to Y$, such that
    \begin{quote}
        \begin{itemize}
            \item[{\rm\bf[CR1]:}] $H(0,\cdot) = f_c$,
            \item[{\rm\bf[CR2]:}] $H(1,x,y) = (L(x),0)$, for all $(x,y)\in N_c$,
            \item[{\rm\bf[CR3]:}] $H(t,x,y) \notin M_c$, for all $(x,y)\in \partial{B_u}\oplus T_N$, $t\in[0,1]$ and
            \item[{\rm\bf[CR4]:}] $H(t,x,y) \in \mathbb{R}_u\oplus T_M$ for $(x,y)\in N_c$ and $t\in[0,1]$.
        \end{itemize}
    \end{quote}
\end{definition}

\begin{figure}
    \centerline{\includegraphics[width=.75\textwidth]{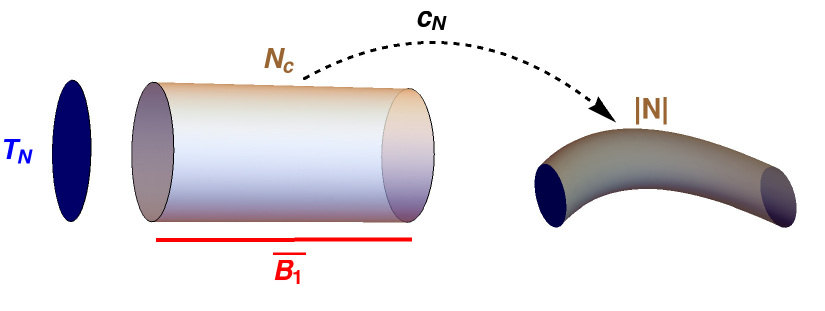}}
    \caption{An example of an h-set in three dimensions with $u(N)=1$ and $T_N=D_2$ -- a two-dimensional closed disc. Here $N_c=\overline{B_1}\oplus D_2$.}
    \label{pic:magicset}
\end{figure}

\begin{figure}
    \centerline{\includegraphics[width=.75\textwidth]{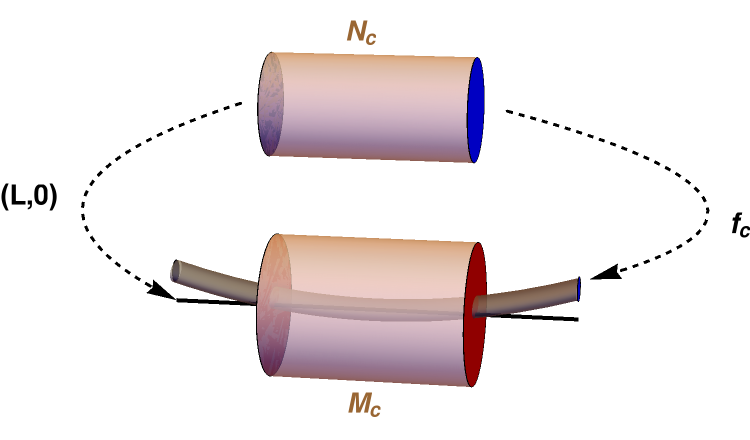}}
    \caption{An example of an $f-$covering relation: $N\cover{f}M$. In this case, the homotopy joining $f_c(x,y)$ with a linear map $(L(x),0)$ and satisfying \textbf{[CR1]}--\textbf{[CR4]} is simply given by $H(t,x,y) = t(L(x),0) + (1-t)f_c(x,y)$.}
    \label{fig:cov}
\end{figure}

A typical picture of a h-set with $u(N)=1$ is given in Figure \ref{pic:magicset}. A  picture illustrating covering relation with one exit direction is given on Figure~\ref{fig:cov}.

The following theorem is the main topological tools used in the proof of the existence of symbolic dynamics in the KS equation \cite[Theorem 1]{WZ}.
\begin{theorem}\cite[Theorem 4, Theorem 6]{WZ}\label{thm:covrel}
    Let $X_i$, $i\in\mathbb Z$ be real normed spaces and let $N_i$ be $h$-sets in $X_i$, respectively. Assume
    \begin{equation*}
    N_i\cover{f_i}N_{i+1}\qquad \text{for }i\in\mathbb Z.
    \end{equation*}
    Then there exists $(u_i)_{i\in\mathbb Z}$ such that
    \begin{eqnarray*}
    u_i\in |N_i|\quad \text{for } i\in\mathbb Z,\\
    f(u_i)=u_{i+1}\quad \text{for } i\in\mathbb Z.
    \end{eqnarray*}
    Moreover, if the sequence of covering relations is periodic, that is there is $K>0$ such that $N_{i}=N_{i+K}$ and $f_{i}=f_{i+K}$ then the sequence $(u_i)_{i\in \mathbb Z}$ can be chosen as $K$-periodic, too.
\end{theorem}

\subsection{Cone condition}
%\label{subsec:cones}
The aim of this section is to describe a technique we use to establish a  convergence of a trajectory of a discrete dynamical system to a periodic trajectory.
For this purpose we adopt the notion of cone condition from \cite{KWZ,ZCC}  to context  of real normed spaces.

We start from two definitions from \cite{KWZ,ZCC} adopted to an in\-fini\-te-di\-men\-sional case.
\begin{definition}
    Let $N \subset X$ be an h-set and let $Q:X \to \mathbb{R}$ be homogenous of degree two and continuous. The pair $(N,Q)$ will be called an \emph{h-set with cones}.
\end{definition}

The function $Q$ can be used to define cone fields by the condition $C^+_x =\{y\in X : Q(y-x)>0\}$ and $C^-_x =\{y\in X : Q(y-x)<0\}$.

\begin{definition}\cite[Definition~13]{KWZ}%\label{def:conec}
    Assume that $(N,Q_N),(M,Q_M)$ are h-sets with cones. Assume that $N \cover{f} M$. We say that $f$ satisfies the cone condition with respect to the pair $(N,M)$, if for any $x_1,x_2 \in N_c$ with
    $x_1 \neq x_2$ there holds
    \begin{equation}\label{eq:cone-cond-def}
    Q_M(f_c(x_1) - f_c(x_2)) > Q_N(x_1 -x_2).
    \end{equation}
    Whenever it is convenient, we will also say that the cone conditions are satisfied for the covering relation $N \cover{f} M$, if the above condition is satisfied.
\end{definition}
Condition (\ref{eq:cone-cond-def}) means that the cone field $C^+$ is preserved by the action of $f$. Moreover, we have a kind of expansion and contraction in $C^+$ and $C^-$, respectively.

The following theorem  provides a tool for proving existence of heteroclinic orbits between two periodic hyperbolic solutions. We would like to emphasize, that the assumptions reduce to finite number of inequalities, and thus in principle they can be checked by a computer.

\begin{theorem}\cite[Theorem~3.4]{W}\label{thm:heteroToolTheorem}
    Assume that
    \begin{eqnarray}
        N_0 \cover{f_0} N_1 \cover{f_1} N_2 \cover{f_2} \cdots \cover{f_{n-1}}
        N_n=N_0, \label{eq:N-chain}  \\
        N_0\cover{h_0}M_1 \cover{h_1}M_2\cover{h_2}\cdots \cover{h_{m-2}}M_{m-1}\cover{h_{m-1}}K_0  \label{eq:NMK-chain} \\
        K_0 \cover{g_0} K_1 \cover{g_1} K_2 \cover{g_2} \cdots \cover{g_{k-1}} K_k=K_0, \label{eq:K-chain}
    \end{eqnarray}
    where
    \begin{itemize}
        \item $(N_i,Q_{N_i})$ are h-sets with cones and $f_i$ satisfies the cone condition with
        respect to $(N_i,N_{i+1})$, $i=0,\ldots,n-1$
        \item $(K_i,Q_{K_i})$ are h-sets with cones and $g_i$ satisfies the cone condition with
        respect to $(K_i,K_{i+1})$, $i=0,\ldots,k-1$
    \end{itemize}
    Then there exist
    \begin{enumerate}
        \item a unique periodic orbit $(u_0,\ldots,u_{n-1})$ with respect to
        $(f_0,\ldots,f_{n-1})$ and the sequence $(N_0,N_1,\ldots,N_{n-1})$
        \item a unique periodic orbit $(v_0,\ldots,v_{k-1})$ with respect to
        $(g_0,\ldots,g_{k-1})$ and the sequence
        $(K_0,K_1,\ldots,K_{k-1})$
        \item a full orbit $(x_i)_{i\in\mathbb Z}$ with respect to
        \begin{multline*}
        (\ldots,
        f_0,\ldots,f_{n-1},f_0,\ldots,f_{n-1},\\
        h_0,\ldots,h_{m-1},\\
        g_0,\ldots,g_{k-1},g_0,\ldots,g_{k-1},\ldots)
        \end{multline*}
        and the sequence
        \begin{equation*}
        \left((N_1,N_2,\ldots,N_n)^{\mathbb
            N},M_1,M_2,\ldots,M_{m-1},(K_0,K_1,\ldots,K_{k-1})^{\mathbb
            N}\right)
        \end{equation*}
        such that
        \begin{eqnarray*}
        x_0\in |N_0|\\
        \lim_{t\to\infty} x_{tk+m}= v_0\\
        \lim_{t\to-\infty} x_{tn}= u_0
        \end{eqnarray*}
    \end{enumerate}
\end{theorem}
We omit the proof as the proof from finite-dimensional  case from \cite{W} is good also in the present setting, as only the compactness  and forward iterations are used. 
Instead we just underline main steps of the proof and explain the assumptions. The chain of covering relations (\ref{eq:N-chain}) satisfying cone conditions 
implies the existence of hyperbolic periodic point in $N_0$. The same applies to chain (\ref{eq:K-chain}).   The chain of covering (\ref{eq:NMK-chain}) gives us transition
from $N_0$ to $K_0$.  The cone conditions implies that two orbits that are in the same h-set forward (backward) in time are approach one another.

%% file: applications.tex
\section{Rigorous results about the KS equation}
\label{sec:conn-orb-ks}

In this section we give computer-assisted proofs of Theorem~\ref{thm:connectingOrbits} and Theorem~\ref{thm:attrPerOrbit}. Both proofs use $\mathcal C^1$ algorithm described in Section~\ref{sec:c1algorithm} to obtain bounds on the derivatives of certain Poincar\'e maps.

\subsection{The KS equation in the Fourier basis}%\label{sec:KSeq}
Consider the equation (\ref{eq:KS}) with periodic and odd boundary conditions (\ref{eq:KSbc}) and assume that $u$ is a classical solution to (\ref{eq:KS})--(\ref{eq:KSbc}) given as a convergent Fourier series
\begin{equation}
u(t,x)=\sum_{k=1}^\infty -2a_k(t) \sin(kx).  \label{eq:ks-func-rep}
\end{equation}

The particular form of representation of $u$ given in (\ref{eq:ks-func-rep})  comes from imposing on $u(t,x)=\sum_k u_k(t) e^{ikx}$ periodic and odd boundary conditions. Then we have  $a_k(t)=\mathrm{Im}\, u_k(t) $.

It is easy to see \cite{CCP,ZM} that for sufficiently regular functions  the system (\ref{eq:KS})--(\ref{eq:KSbc})  gives rise to an infinite ladder of coupled ODEs
\begin{equation} \label{eq:fuKS}
\frac{d a_k}{dt}=F_k(a):= k^2(1-\nu k^2) a_k - k \sum_{n=1}^{k-1} a_n
a_{k-n} + 2k \sum_{n=1}^{\infty} a_n  a_{n+k}, \quad k\geq 1.
\end{equation}

Define our space $H=c_0$ - the space of sequence converging to zero with the norm $\|a\|_\infty=\sup_{k}|a_k|$.
On different Poincar\'e sections for the system (\ref{eq:fuKS}) we will use an equivalent norm $\|a\|=\max \left(\sqrt{\sum_{i\leq m} |a_k|^2},\sup_{k >m} |a_k|\right)$ for some $m>0$.

In order to apply Theorem~\ref{thm:covrel} and Theorem~\ref{thm:heteroToolTheorem} discussed in Section~\ref{sec:geometricTools} to a Poincar\'e map for (\ref{eq:fuKS}) we need to chose a family of compact sets. The motivation for choosing sequences with geometrically decaying coefficients is given in \cite{WZ}. Here we introduce the notation.

For $S>0$ and $q>1$ we set
\begin{eqnarray}\label{eq:geometricSetForm}
    \WqS{q}{S} &=& \left\{(a_k)_{k=1}^\infty \,|\, a_k \in \mathbb{R}, \ |a_k|\leq Sq^{-k}\text{ for }k\geq 1\right\}.
\end{eqnarray}
The set $\WqS{q}{S} \subset H$  is compact for any $S>0$ and $q>1$ -- see \cite[Theorem 26]{WZ24}. It is also compact as a subset $\WqS{q}{S} \subset H^s([-\pi,\pi])$ for any $s$ and in the space of analytic $2\pi$-periodic functions on the strip $\mathbb{R} \times [-\delta,\delta]$ for $\delta < q$. In fact a coordinate-wise convergence in $\WqS{q}{S}$  is equivalent to convergence in $C^0$, $H^s$ and in the space of analytic functions. Indeed the equivalent norm in $H^s$ in terms of the Fourier coefficient can be defined as
\begin{equation*}
  \|(a_k)_{k=1}^\infty \|_{H^s}=\left(\sum_{k\geq 1} |k|^{2s} |a_k|^2 \right)^{1/2},
\end{equation*}
and the convergence in such norm on $\WqS{q}{S}$ is equivalent the convergence for each $k$, because the tails of elements from  $\WqS{q}{S}$
are uniformly bounded in $H^s$.

\subsection{Computer--assisted proof of Theorem~\ref{thm:connectingOrbits}.}

In what follows, by a Poincar\'e section $\Theta$ we mean a subset of a codimension one affine subspace $\widetilde \Theta$, such that the vector field is transverse to $\widetilde \Theta$ on $\Theta$. In the applications the sets $\Theta$ will be intersections of $\widetilde \Theta$ with the sets of the form $W_{q,S}$. Let $\Theta_1,\Theta_2$ be Poincar\'e sections for (\ref{eq:fuKS}). By $\PM_{\Theta_1 \to \Theta_2}$ we denote the Poincar\'e map between two sections, that is for $a\in\Theta_1$ the point $\PM_{\Theta_1 \to \Theta_2}(a)$ is the first intersection of the trajectory of $a$ with $\Theta_2$, provided it exists. We also define
$$\PM_{\Theta_1\to\Theta_2}^i= \PM^{i-1}_{\Theta_2 \to \Theta_2}  \circ \PM_{\Theta_1 \to \Theta_2}.  $$

The following lemma has been proven with computer assistance and it was the crucial step in the proof of Theorem~\ref{thm:symdynKS}.

\begin{lemma}\cite[Lemma~10]{WZ}\label{lem:symbolic-dynamics}
    Consider the system (\ref{eq:fuKS}) on $\ltwoq$ with $q=\frac{3}{2}$. There are explicitly given Poincar\'e sections $\Theta^1$, $\Theta^2$ and $\Theta^{i}_{1\to2}$, $\Theta^i_{2\to 1}$,   $i=0,\ldots,10$ and explicitly given pairwise disjoints $h$-sets $N^j\subset \Theta^j$, $j=1,2$, $N^i_{1\to 2}\subset \Theta^{i}_{1\to 2}$, $N^i_{2\to 1}\subset \Theta^{i}_{2\to 1}$, $i=0,\ldots,10$ such that
    \begin{multline*}
    N^1\cover{\PM^2}N^1 \cover{\PM^3} N^0_{1\to 2} \cover{\PM^5} N^1_{1\to 2}\cover{\PM^5}N^2_{1\to 2}\cover{\PM^5}\cdots \cover{\PM^5} N^{10}_{1\to 2}\cover{\PM^5} N^2,\\
    N^2\cover{\PM^4}N^2 \cover{\PM^4} N^0_{2\to 1} \cover{\PM^5} N^1_{2\to 1}\cover{\PM^4}N^2_{2\to 1}\cover{\PM^5}\cdots \\
    \cdots \cover{\PM^4} N^{8}_{2\to 1}\cover{\PM^6} N^{9}_{2\to 1}\cover{\PM^5} N^{10}_{2\to 1}\cover{\PM^3} N^1,
    \end{multline*}
    where the starting and target sections in each of the above covering relations are determined by the h-sets appearing in the relation.
\end{lemma}

\begin{figure}[htp]
    \centerline{\includegraphics[width=.9\textwidth]{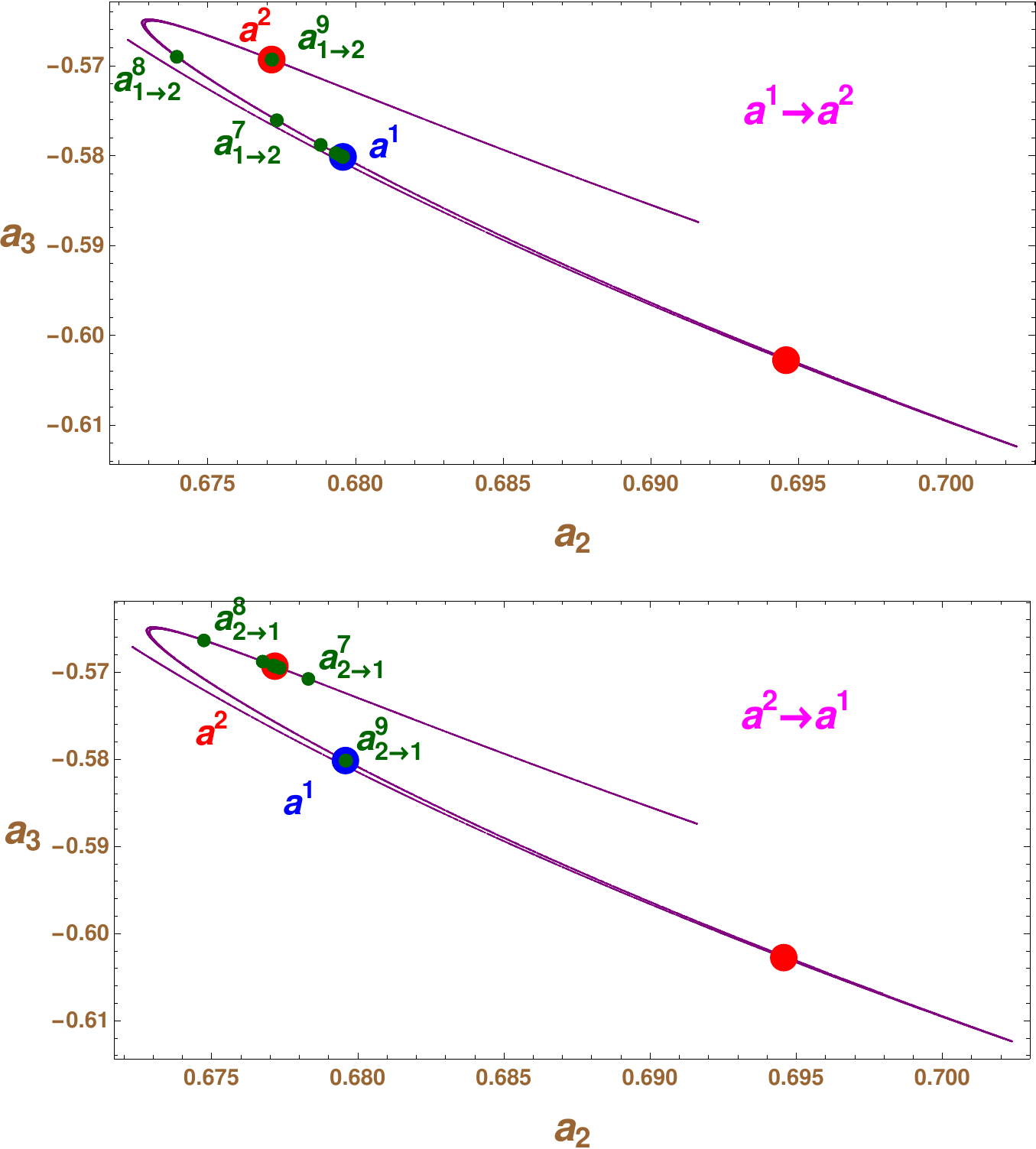}}
    \caption{Numerically observed heteroclinic connections $a^{i}_{1\to 2}$, $a^{i}_{2\to 1}$ $i=0,\ldots,10$ between approximate periodic points $a^1$ and $a^2$. All these points are located in the section $\Theta=\{a\in l_2 : a_1=0 \wedge a_1'=f_1(a)>0\}$ .\label{fig:heteroclinic}}
\end{figure}

The sets $N^j$, $j=1,2$, are centred at approximate periodic points $a^j$ of periods $j$, respectively, which correspond to periodic solutions $u^j$ of the original system (\ref{eq:KS}) with boundary conditions (\ref{eq:KSbc}) --- see Figure~\ref{fig:periodic-orbits-3d}. The sets $N^i_{j\to k}$ are centred at approximate heteroclinic points $a^{i}_{1\to2}$, $a^{i}_{2\to1}$, $i=0,\ldots,10$, found by extensive numerical studies --- see Figure~\ref{fig:heteroclinic}.

In order to obtain connecting orbits between periodic points in $N^1$ and $N^2$ it suffices to check cone conditions for $N^1\cover{\PM^2}N^1$ and $N^2\cover{\PM^4}N^2$ and apply Theorem~\ref{thm:connectingOrbits}. Preliminary simulation showed however, that in order to obtain bounds on derivatives sharp enough to check the cone conditions on these sets we would need to subdivide them into very large number of smaller sets (note we have to subdivide in infinite-dimensional space). This would cause very large time of computation.

\begin{figure}[htbp]
    \centerline{\includegraphics[width=\textwidth]{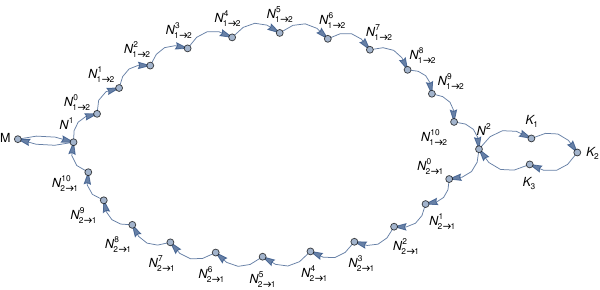}}
    \caption{Graph of symbolic dynamics for the KS equation.\label{fig:graphsymbolic-dynamics}
        }
\end{figure}

To avoid this problem, we introduced intermediate sections and constructed loops of covering relations satisfying cone condition, as shown in Figure~\ref{fig:graphsymbolic-dynamics}. The centres of new sets in these loops, that is $M, K_1,K_2$ and $K_3$ are located at intersections of corresponding approximate periodic orbits with the Poincar\'e section $a_1=0$.

Proof of the following lemma is computer-assisted.
\begin{lemma}\label{lem:KS-cone-conditions}
    There exist explicitly given Poincar\'e sections $\Theta^1_1$, $\Theta^2_1$, $\Theta_2^2$ and $\Theta_3^2$
    \begin{itemize}
        \item $h$-set with cones $N^1\subset\Theta^1$ and $h$-set with cones $M\subset \Theta_1^1$,
        \item $h$-set with cones $N^2\subset \Theta^2$ and $h$-sets with cones $K_i\subset \Theta^2_{i}$, $i=1,2,3$,
    \end{itemize}
 such that
 \begin{eqnarray}
 M\cover{\PM}N^1\cover{\PM}M,\label{eq:cover-1}\\
 N^2\cover{\PM}K_1\cover{\PM}K_2\cover{\PM}K_3\cover{\PM}N^2.\label{eq:cover-2}
 \end{eqnarray}
Moreover, the cone conditions are satisfied for these covering relations.
\end{lemma}

\textbf{Proof:}
    The trajectory of $a^1$, which is an approximate fixed point for $\PM^2$, intersects the section $a_1=0$ at $a^1_1\approx \PM(a^1)$. At $a^1_1$ we attached a Poincar\'e section as a hyperplane almost orthogonal to the vector field at $a_1^1$. Then we defined coordinate system in $\Theta_1^1$ as approximate eigenvectors of $D\PM^2(a_1^1)$. In this coordinate system we defined the set $M$ as a parallelogram centred at $a_1^1$ and with edges parallel to the axes.

    Similarly, we set $a^2_i \approx \PM^i(a^2)$, $i=1,2,3$. At each of these points we attached a Poincar\'e section $\Theta_i^2$ as a hyperplane almost orthogonal to the vector field at $a_i^2$. Then we defined coordinate systems in $\Theta^2_i$ as approximate eigenvectors of $D\PM^4(a_i^2)$, $i=1,2,3$, respectively. Finally, we defined the sets $K_i\subset \Theta^2_i$ as parallelograms centred at $a^2_i$.

    Then, using algorithm for verification of covering relations \cite[Appendix A.4]{WZ} we checked (\ref{eq:cover-1})--(\ref{eq:cover-2}).

    For verification of cone conditions we set homogeneous degree two continuous functions associated with $N^1$, $N^2$, $M$, $K_1$, $K_2$ and $K_3$ as
    \begin{equation*}
    \begin{array}{rcl}
        Q_{N^1}(a) &=& a_1^2 - 0.5a_2^2 - a_3^2-\cdots - a_{17}^2  - \|y\|_\infty^2,\\
        Q_{M}(a) &=& 2a_1^2 - a_2^2 - a_3^2-\cdots - a_{17}^2  - \|y\|_\infty^2,\\
        Q_{N^2}(a) &=& a_1^2 - a_2^2 - a_3^2-\cdots - a_{17}^2  - \|y\|_\infty^2,\\
        Q_{K_1}(a) &=& a_1^2 - a_2^2 - a_3^2-\cdots - a_{17}^2  - \|y\|_\infty^2,\\
        Q_{K_2}(a) &=& a_1^2 - a_2^2 - a_3^2-\cdots - a_{17}^2  - \|y\|_\infty^2,\\
        Q_{K_3}(a) &=& 4a_1^2 - 2a_2^2 - a_3^2-\cdots - a_{17}^2  - \|y\|_\infty^2,\\
    \end{array}
    \end{equation*}
    Note, in each case the arguments of the above functions are given in section coordinates specific for the $h$-set. To validate cone conditions we use Lemma~\ref{lem:rigConeVerify} from the Appendix. To this end, we computed quantities $a$, $d$ and $c$ as defined in Lemma~\ref{lem:rigConeVerify}. Obtained bounds are listed in Table~\ref{tab:eigenvalues}.
\begin{table}
    \begin{center}
    \begin{tabular}{c|c|c|c|c}
        covering relation & $a$ & $d$ & $c$ & $ad -c^2$ \\ \hline
        $N^1\cover{\PM}M$ & $\geq 0.009$ & $\geq 0.999$ & $\leq 0.025$ & $>0.008$\\ \hline
        $M\cover{\PM}N^1$ & $\geq 0.469$ & $\geq 0.997$ &  $\leq 0.042$ & $>0.4$\\ \hline
        $N^2\cover{\PM}K_1$ & $\geq 0.626$ & $\geq 0.999$ & $\leq 0.008$ & $>0.6$\\ \hline
        $K_1\cover{\PM}K_2$ & $\geq 0.882$ & $\geq 0.991$ & $\leq 0.103$ & $>0.8$\\ \hline
        $K_2\cover{\PM}K_3$ & $\geq 0.041$ & $\geq 0.997$ & $\leq 0.1$ & $>0.03$\\ \hline
        $K_3\cover{\PM}N^2$ & $\geq 0.142$ & $\geq 0.995$ & $\leq 0.063$ & $>0.1$\\ \hline
    \end{tabular}
    \caption{Bounds on $a$, $d$ and $c$ defined as in Lemma~\ref{lem:rigConeVerify} for all covering relations in Lemma~\ref{lem:KS-cone-conditions}. \label{tab:eigenvalues}}
    \end{center}
\end{table}
    From these estimates it is clear, that in each case $a>0$ and $ad>c^2$ and thus the assumptions of Lemma~\ref{lem:rigConeVerify} and in consequence the cone condition are satisfied.

Summarizing, the assumptions of Theorem~\ref{thm:heteroToolTheorem} are satisfied for the set of covering relations from Lemma~\ref{lem:symbolic-dynamics} and (\ref{eq:cover-1})--(\ref{eq:cover-2}). This completes a computer-assisted proof of Theorem~\ref{thm:connectingOrbits}.
\qed

\subsubsection{Implementation notes}

The computation related to the existence of all covering relations from \\ Lemma~\ref{lem:symbolic-dynamics}
and (\ref{eq:cover-1})--(\ref{eq:cover-2}) except $N^i\cover{\PM^{2i}}N^i$, $i=1,2$, took approximately  95 minutes on a standard personal computer with AMD Ryzen 7 5800X 3.80GHz processor and running 16 concurrent threads. This is slightly faster than the time reported in \cite{WZ}, that is 40 minutes on 78 CPUs.

We observed, that the strategy of choice of one of the parameters in the automatic differentiation procedure that computes Taylor coefficients of the solutions to (\ref{eq:fuKS}) has a  huge impact on obtained bounds. This allowed us to enlarge time steps of the solver and in consequence shorten the time of computation.

The vector $F_k(a)$ when evaluated on a set $a$ with geometric decaying tail $|a_k|\leq Cq^{-k}$ for $k>m$ and some $C\geq 0$, $q>1$ can be bounded by a polylogarithm function, that is $|F_k(a)|\leq \widetilde C k^4q^{-k}$ for some computable $\widetilde C\geq0$. In \cite[Lemma 28]{WZ} we shown that if we fix any $\delta\in(1,q)$, then $|F_k(a)|\leq D (q/\delta)^{-k}$ for a computable $D\geq 0$ and thus we obtained desired geometric decay of $F_k$. It turned out that the obtained bounds are highly sensitive to the choice of $\delta$. The current improvement has been done by trial and choice of the strategy for $\delta$.

This problem disappears when we use instead a larger family of sequences to represent $a$ and $F(a)$, that is $|a_k|\leq C k^sq^{-k}$ for some $C\geq 0$, $q>1$ and $s\in\mathbb Z$. On the other hand, algebraic operations on such objects will be slower. Although we did not made detailed analysis in this direction we believe, that overall this is the right way to further reduce overestimation and speed up computation.

All above comments apply to computation of variational equations. First version of the program that verified the cone condition run over $10$ hours on a computer with $448$ CPUs. After these changes, and most important introduction of intermediate sections (that is sets $K$, $M_1,M_2,M_3$) we reduced the time of computation related to cone condition to approximately an hour on a personal computer running 12 concurrent threads.

The program have been written in C++ and is based on the CAPD library \cite{CAPDREVIEW}. The source code can be found at \cite{WSc}.

\subsection{Attracting periodic orbit for $\nu=0.127$}
\label{sec:attr-per-orb}

The goal of this section is to give a computer assisted proof of Theorem~\ref{thm:attrPerOrbit}. The system (\ref{eq:fuKS}) admits a symmetry
$$S(a_1,a_2,a_3,a_4,\ldots) = (-a_1,a_2,-a_3,a_4,\ldots).$$
Let $\Theta=\{(a_i)_{i=1}^\infty : a_1=0\}$ be a Poincar\'e section. Since $\Theta = S(\Theta)$ it follows that $S$ is also a symmetry for the associated Poincar\'e map $\PM$ and the same holds true for Poincar\'e maps for Galerkin projections. Since $S$ is an involution, we have $\PM^2 = (S\circ \PM)^2$. Hence, it suffices to prove the existence of an attracting fixed point for $S\circ \PM$.

The existence of a fixed periodic point for $S\circ \PM$ is checked by means of the Schauder fixed point theorem. For this we need to define a compact and convex set $K\subset \Theta$ and check that $\widetilde K := S\circ \PM(K)\subset K$. Then we will check that $\|D(S\PM)(\widetilde K)\|<1$. Observe that this implies that $K$ is contraction on $K$.

Our space $H=c_0$ - the sequences converging to zero.
Let us fix $m=18$. We define norm on $H$ by
\begin{equation}
 \|u\|=\max \left(\left(\sum_{j=1}^m |u_j|^2\right)^{1,2},\sup_{j>m} |u_j| \right). \label{eq:atr-norm}
\end{equation}
It is immediate that $H$ with above norm is \gss space.

Let $\PM_m$ denote Poincar\'e map for the flow induced by $m$-dimensional Galerkin projection. The point $u^0=(u_1,\ldots,u_m)\in\Theta_m$, where
$$
\begin{array}{l|l}
u^0_{1}=0 & u^0_{2}=1.1655388447652766\\
u^0_{3}=5.7913642084409762\cdot 10^{-1} & u^0_{4}=-2.768914122643753\cdot 10^{-1}\\
u^0_{5}=-1.2582912341558827\cdot 10^{-1} & u^0_{6}=1.3015613792264165\cdot 10^{-2}\\
u^0_{7}=1.6757421753649995\cdot 10^{-2} & u^0_{8}=7.3178705708507339\cdot 10^{-4}\\
u^0_{9}=-1.4755994209479948\cdot 10^{-3} & u^0_{10}=-2.5601313918591901\cdot 10^{-4}\\
u^0_{11}=9.5427473299718919\cdot 10^{-5} & u^0_{12}=3.2627521780891864\cdot 10^{-5}\\
u^0_{13}=-3.7164835203789059\cdot 10^{-6} & u^0_{14}=-2.9884428114398931\cdot 10^{-6}\\
u^0_{15}=-6.6059155560432874\cdot 10^{-8} & u^0_{16}=2.1596485662114442\cdot 10^{-7}\\
u^0_{17}=2.9759304684536878\cdot 10^{-8} & u^0_{18}=-1.2262214633625861\cdot 10^{-8}\\
\end{array}
$$
is an approximate symmetric fixed point for $S\circ \PM_m$ (period two point for $\PM_m$, that is
$$
\|S(\PM_m(u^0)) - u^0\|_2\approx 4\cdot10^{-13}
$$
and it has been found by a simple iteration.

On the Poincar\'e section $\Theta$ we introduce an affine coordinate $u=(u^0 + A\cdot P_mr,Q_mr)$, where $A$ is a matrix that almost diagonalizes $SD\PM_m(u)$ (nonrigorous computation) and $P_m,Q_m$ are projections. In what follows we will work in this coordinate system, that is $r=0$ corresponds to $(u^0,0)$.

As a preliminary step we compute rigorously (using $\mathcal C^0$ algorithm from \cite{WZ}) a bound $K^0\supset S(\PM(0))$. Clearly any reasonable choice of the set $K$ should contain $K^0$. As a candidate for $K$ we take
$$K_i = K_i^0\cdot [-10,10],\quad i=1,\ldots, m$$
(note $K_1^0=0$) and the tail is blown up by a factor $1.25$, that is $|K_i|\leq 1.25Cq^{-i}$, where $C$ is the constant bounding tail in $K^0$, that is $|K^0_i|\leq Cq^{-i}$ for $i>m$. Using $\mathcal C^0$ algorithm form \cite{WZ} we compute a bound  $\widetilde K\supset S\PM(K)$ (in the coordinate system mentioned above) and we get
\begin{equation}\label{eq:fixPointC0Bound}
\begin{array}{|c|c|c|}
    i & \max |K_i|  & \max |\widetilde K_i| \\ \hline
2 & 1.4037896821825919\cdot 10^{-8} & 8.7024122740689425\cdot 10^{-9}\\
3 & 4.9756335159569671\cdot 10^{-9} & 8.5181999007451079\cdot 10^{-10}\\
4 & 9.7156048870811418\cdot 10^{-11} & 6.6221478045718261\cdot 10^{-12}\\
5 & 2.637578699436423\cdot 10^{-11} & 2.6821395044388368\cdot 10^{-12}\\
6 & 8.491763367654622\cdot 10^{-12} & 8.6150267406921159\cdot 10^{-13}\\
7 & 1.6209035540532317\cdot 10^{-11} & 1.7003223001285466\cdot 10^{-12}\\
8 & 1.8203180667485086\cdot 10^{-11} & 1.8342179335531725\cdot 10^{-12}\\
9 & 4.7853573165474597\cdot 10^{-11} & 4.8670536302275007\cdot 10^{-12}\\
10 & 2.677612693142676\cdot 10^{-10} & 2.6635853628312471\cdot 10^{-11}\\
11 & 5.6520226737222243\cdot 10^{-11} & 5.6692000787577923\cdot 10^{-12}\\
12 & 1.0201585471030834\cdot 10^{-10} & 1.013058647787179\cdot 10^{-11}\\
13 & 4.341777521371891\cdot 10^{-10} & 4.3005552049615385\cdot 10^{-11}\\
14 & 1.4398858307285336\cdot 10^{-10} & 1.4274077590841134\cdot 10^{-11}\\
15 & 3.3157949280900669\cdot 10^{-10} & 3.2774891340532832\cdot 10^{-11}\\
16 & 7.7478332603484965\cdot 10^{-10} & 7.6280987034794698\cdot 10^{-11}\\
17 & 1.0424399029461924\cdot 10^{-9} & 1.0253441358095485\cdot 10^{-10}\\
18 & 8.5396225348412251\cdot 10^{-10} & 8.3748052013014238\cdot 10^{-11}\\ \hline
C & 9.9201773737522156\cdot 10^{-6} & 7.9361550347737575\cdot 10^{-6}\\

\end{array}.
\end{equation}

Because $K_i$ is symmetric around zero it is clear from the above data that $\widetilde K\subset K$ and thus the existence of a fixed point $u^*\in \widetilde K$ for $S\circ \PM$ is proved.

In the second step we compute a bound on the operator norm of $A:=SD\PM(\widetilde K)$. On Poincar\'e section $\Theta= X\times Y$ we use the norm $\|(x,y)\| = \max\left\{\|x\|_2,\|y\|_\infty\right\}$, where $X$ is $m-1$-dimensional (the section is given by $a_1=0$) and $Y$ is the tail. Using the $\mathcal C^1$ algorithm from Section~\ref{sec:c1algorithm} we computed bounds
\begin{equation}\label{eq:fixPointC1Bound}
    \begin{array}{rcl}
    \|A_{xx}\| &\in&  [0.52569744154488263, 0.52578830885349082]\\
    \|A_{xy}\| &\leq& 0.0021336347518391241 \\
    \|A_{yx}\| &\in& [5.6323483625902106\cdot 10^{-321}, 7.8313759074032856\cdot 10^{-18}]\\
    \|A_{yy}\| &\leq&7.8313759074032856\cdot 10^{-12}.
\end{array}
\end{equation}
From the above data and the definition of the norm (see (\ref{eq:atr-norm}))  it is clear, that $$\|A\| = \max\left\{ \|A_{xx}\|+\|A_{xy}\|,\|A_{yx}\|+\|A_{yy}\|\right\} \leq 0.52792194360532996 <1.$$
\qed
\begin{remark}
The computation regarding data in (\ref{eq:fixPointC0Bound}) and (\ref{eq:fixPointC1Bound}) took $11$ minutes and $73$ minutes, respectively, in a single thread on laptop-type computer with i7-8650U, 1.90GHz processor. The computation were performed in effective dimension $m=18$ and with the tolerance (parameter of the algorithm specifying the acceptable error per one time step of integration) set to $\texttt{\emph{tol}}=10^{-11}$.

Changing the setting to $m=15$ and $\texttt{\emph{tol}}=10^{-7}$ the validation also succeeded within $26$ seconds and $217$ seconds for the $\mathcal C^0$ and $\mathcal C^1$ parts, respectively. Obtained bound on the norm is much wider, that is
$$\|A\|\leq 0.92265954650065929,$$
although the norms in both cases are different ($m=18$ and $m=15$ give different decomposition of the space). On the other hand, the stability domain is much larger, as the set $K$ on which we perform validation is much bigger ($4$ orders of magnitude).
\end{remark}

%\textbf{Szkic:}
%\begin{itemize}
% \item the existence of compact absorbing set $\mathcal{K}$
% \item smoothing of solutions in $\mathcal{K}$, after finite time we enter $\mathcal{A'}$ which correct decay - of our type
% \item from the continuity of the evolution of KS  somehow we should be able to find $\mathcal{U}$ and $\mathcal{A}$ as stated in the assertion
%\end{itemize} 

%% file: varsysintro.tex
\section{Variational equations in infinite dimensions}
\label{sec:variationalEquations}

The goal of this section present some results from~\cite{WZ24} about Galerkin projections of variational equations for class of dissipative PDEs. This is needed to
justify our $C^1$-algorithm for KS equation discussed in Section~\ref{sec:c1algorithm}.

By a $\mathcal C^1$ algorithm for an ODE we understand computation of bounds on the solutions to the following initial value problem
\begin{alignat}{3}
	&\frac{d}{dt}x(t)=F(x(t)),\qquad &&x(0)=x_0\in H, \label{eq:c0IVP}\\
	&\frac{d}{dt}V(t,x)=DF(x(t))V(t,x), \qquad &&V(0,x_0)=V_0\in \mathrm{Lin}(H,H),\label{eq:c1IVP}
\end{alignat}
where $H$ is a real Banach space and $F\colon H\to H$. It is well known that if $H=\mathbb R^n$ and $F\in\mathcal C^2(\mathbb R^n)$, the IVP (\ref{eq:c0IVP})-(\ref{eq:c1IVP}) admits unique solution. Moreover, if $V_0=\mathrm{Id}_n$ then $V(t,x_0)=\frac{d}{dx}\varphi(t,x_0)$ is the partial derivative with respect to initial condition of the induced local flow and we have the following useful Euler formula
\begin{equation}\label{eq:MVFforFlows}
	\varphi(t,x)-\varphi(t,y) = \int_0^1 V(t,y+s(x-y))ds\cdot (x-y),
\end{equation}	
provided the right hand side is defined.

Extension of (\ref{eq:c0IVP})-(\ref{eq:MVFforFlows}) to infinite dimension is, however, more delicate. In many interesting applications the vector field $F:H\to H$ is defined on a subspace of a real Banach space $H$, only. Also the Fr\'echet derivatives of the vector field and/or induced semi-flow may not be defined. On the other hand, (\ref{eq:c0IVP})-(\ref{eq:MVFforFlows}) is well posed for any $n$-dimensional $\mathcal C^2$-smooth approximation (Galerkin projection) of the original infinite-dimensional problem. Hence, there is a chance to obtain a formula similar to (\ref{eq:MVFforFlows}) by passing to the limit with these approximations. %, provided the limit operators $V(t,y+s(x-y))\in\mathrm{Lin}(H,H)$ are properly defined and understood.

This question was discussed by the authors in  \cite{WZ24} and in the rest of this section we recall from there the relevant definitions and theorems.

%\textbf{Before we describe our algorithm for rigorous integration of variational equations for KS equation we need to } recall  from \cite{WZ24} the most important results about convergence of solutions to finite-dimensional Galerking projections and associated variational equations to the solutions of underlying infinite-dimensional. These abstract results (\cite[Thm. 11 and Thm. 14]{WZ24}) provide a framework for the algorithm for integration of variational equations for a class of dissipative PDEs. From these results we conclude, that the objects we would like to enclose indeed exists inside sets constructed by the algorithm. Informally speaking, if we can enclose uniformly solutions to all finite-dimensional Galerkin projections, then these uniform bounds contain solutions to underlying infinite-dimensional PDE. The algorithm itself will be presented in Section~\ref{sec:c1algorithm}.

\subsection{Notation}
By $\mathbb{Z}_+$ we denote the set of all positive integers. Given a matrix $A \in \mathbb{R}^{n \times n}$ and a norm $\|\cdot\|$ by
\begin{equation*}
	\mu(A)=\lim_{h \to 0^+} \frac{\|\mathrm{Id} + hA \| - 1}{h}
\end{equation*}
we denote the logarithmic norm of $A$ induced by the norm $\|\cdot\|$ --- see \cite{D58,L58,KZ} for the properties. In particular, for a frequently used in this article maximum norm  $\|x\|_{\infty}=\max_{i=1,\dots,n}  |x_i|$ on $\mathbb{R}^n$, there holds
\begin{equation*}
	\mu(A)=\max_{i=1,\dots,n} \left( A_{ii} + \sum_{k=1 , k\neq i}^n |A_{ik}| \right). %\label{eq:logNormMax}
\end{equation*}

For a (possibly infinite dimensional) matrix $V$ by $V_{i\ast}$ and $V_{\ast j}$ we denote its $i$-th row and $j$-th column, respectively.

Given two normed vector spaces $V,W$ by $\text{Lin}(V,W)$ we will denote the space of all bounded linear maps from $V$ to $W$.

\subsection{Convergence conditions}
%\label{sec:convergenve-conditions}

Let us fix a real Banach space $H \subset \mathbb{R}^{\mathbb{Z}_+}$, i.e. it is a space of sequences with some norm $\|\cdot\|$. Let $e_i$ denote the sequence such that $(e_i)_j=0$ for $j \neq i$ and $(e_i)_i=1$. For a nonempty set $J \subset \mathbb{Z}_+$ by $P_J$ we denote a projection defined by $P_J(w)=\sum_{i\in J}w_i e_i$ and by $P_n$ we denote the projection $P_n(\{w\}_{i\in\mathbb Z_+})=w_1e_1+\ldots +w_ne_n$.

For $n \in \mathbb{Z}_+$  by $H_n$ we denote a subspace spanned by
$\{e_j\}_{j \leq n}$. Put $P_n:=P_{\{j \leq n\}}$ and $Q_n=\mathrm{Id}-P_n$. By $\iota_n:H_n \to H$ we denote the embedding $H_n$ into $H$.

We define our standing assumptions \NC1--\NC5 on the space $(H,\|\cdot\|)$.
\begin{definition}\cite[Def. 1]{WZ24}
	%\label{def:gss-norm}
	We will say that $H \subset \mathbb{R}^{\mathbb{Z}_+}$ with norm $\| \cdot \|$ is \emph{\gss} (good sequence space) if the following conditions are satisfied.
	\begin{description}
		\item[\NC 1] $H$ is a Banach space.
		\item[\NC 2] $\forall w \in H\ w=\sum_i w_i e_i$.
		\item[\NC 3] For all  $w \in H$ and for any $\alpha \in \{-1,1\}^{\mathbb{Z}_+}$ there holds $w^\alpha=\sum_i \alpha_i w_i e_i \in H$
		and $\|w\|=\|w^\alpha\|$.
		\item[\NC 4] $\|P_J w\| \leq \|w\|,  \quad \forall w \in H, \, \forall J \subset \mathbb{Z}_+$.
		\item[\NC 5]  there exists  constant $G$, such that for all $w \in H$ and for all $i$ $|w_i| \leq G \|w\|$.
	\end{description}
\end{definition}

\textbf{Examples:} $l_2$, $l_1$, $c_0$(sequences converging to $0$)  with the norm $\|\cdot\|_\infty$. However $l_\infty$ is not \gss space, because \NC{2} is not satisfied. The space defined by convergence of $\sum_i |w_i|2^i$ does not satisfy \NC{5}.

By $H_n$ we denote the subspace of $H$ spanned by $\{e_i\}_{i=1,\ldots,n}$.  The $n$-th Galerkin projection for  (\ref{eq:c0IVP})-(\ref{eq:c1IVP}) is defined as
\begin{equation}\label{eq:varSys}
	\begin{cases}
		\frac{d}{dt}x^n(t)=P_nF(x^n(t)),\\ \frac{d}{dt}V^n(t,x^n)=DP_nF_{|H_n}(x^n(t))V^n(t,x^n), \\ x^n(0)=P_nx_0\in H_n, \\ V^n(0,P_nx_0)=P_nV_0\in \mathrm{Lin}(H_n,H_n).
	\end{cases}
\end{equation}
%We assume that the components $F_i$ are formal power series in variables $a_k$ and symbols $\frac{\partial^k F_i}{\partial a_{i_1} \dots \partial a_{i_k}}$ are defined by formal differentiation of $F_i$.

\begin{definition}\cite[Def. 4]{WZ24}
	Assume that $(H,\|\cdot\|)$ is \gss space. We say that $W\subset H$ satisfies condition $\mathbf{S}$ if
	\begin{description}
		\item[S1:] $W$ is convex and there exists $M \geq 1$, such that $P_n(W) \subset W$
		for $n \geq M$,
		\item[S2:] $W$ is compact.
	\end{description}
\end{definition}

Define $F^n:H_n \to  H_n$ by $F^n(z)=P_n (Fz)$, i.e. this is simply a Galerkin projection of the vector field $F$.
\begin{definition}\cite[Def. 5]{WZ24}
	We say that $F$ is admissible, if $F^n\in\mathcal{C}^3(H_n)$ for all $n>0$.
\end{definition}

\begin{definition}\cite[Def. 6]{WZ24}
	%\label{def:C1C2C3}
	Let $(H,\|\cdot\|)$ be \gss,  $W\subset H$ and $F:\mathrm{dom}(F)\subset H \to H$. We say that $F$ satisfies convergence condition \textbf{C} on $W$ if $F$ is admissible, $W$ satisfies condition \textbf{S} and
	\begin{description}
		\item[C1:] $W \subset \dom(F)$ and function $F|_W:W \to H$ is continuous;
		\item[C2:] there exists  $l \in \mathbb{R}$ such that for all $n\in\mathbb Z_+$ there holds
		\begin{equation*}
			\sup_{x \in P_n W}\mu\left(D F^n (x) \right) \leq l.
		\end{equation*}
	\end{description}
\end{definition}

Let $f:H \supset \dom (f) \to \mathbb{R}$. In the applications we keep in mind, the function $f$ will be a component of a vector field, which is often Frech\'et non-differentiable. On the other hand, assuming that the restrictions $f^n:=f\circ i_n: H_n \to \mathbb{R}$ are $\mathcal C^1$-smooth for all $n$, the derivative $Df^n(z)$ can be seen as a linear form on $H_n$ represented by a row-vector $\left(\frac{\partial f}{\partial x_1}(z),\dots,  \frac{\partial f}{\partial x_n}(z)\right)$. Formally $Df^n(z) \in H_n^*$ for $z \in H_n$. We can embed it to $H^*$ by setting $Df^n(z) \circ P_n $ for $z \in H_n$.

Here a natural question arises: does the infinite row-vector \newline
$\left(\frac{\partial f}{\partial x_1}(z), \frac{\partial f}{\partial x_2}(z),\dots\right)$ define a $1$-form on $H$? If so, can we treat it as the derivative of $f$? For this purpose we introduce the notion
\begin{equation*}%\label{eq:rowderivative}
	\widetilde{D f}(z) := \lim_{n \to \infty} Df^n(P_n z) \circ P_n.
\end{equation*}

Now we can define a vector field on $H \times H$ by
\begin{equation}
	\frac{d}{dt} (x,C)=F_V(x,C):=(F(x),\{\widetilde{DF_i}(x)C\}_{i\in\mathbb Z_+}).  \label{eq:sysVarH}
\end{equation}

\begin{definition}\cite[Def. 7]{WZ24} %\label{def:condV}
	Let $(H,\|\cdot\|)$ be \gss,  $W\times W_V\subset H\times H$ and $F:\mathrm{dom}(F)\subset H \to H$. We say that $F$ satisfies convergence condition \textbf{VC} on $W\times W_V$ if $F$ is admissible, $W\times W_V$ satisfies condition \textbf{S} and
	\begin{itemize}
		\item[\bf{VC1:}] the function $F_V$ defined in (\ref{eq:sysVarH}) is continuous on $W\times W_V$;
		\item[\bf{VC2:}] there exists a constant A, such that
		$$
		\sup_{(x,v)\in P_n(W\times W_V)}\mu(DP_nF_V|_{H_n\times H_n})\leq A.
		$$
	\end{itemize}
\end{definition}
The following theorem  summarizes the most important consequences of the convergence conditions \textbf{C} and \textbf{VC}.
\begin{theorem}{\cite[Thm. 11 and Thm. 14]{WZ24}}
	\label{thm:limitLN}
	Assume that:
	\begin{enumerate}
	\item $(H,\|\cdot\|)$ is \gss
	\item $F$ satisfies {\rm \textbf{C}} on $W\subset H$
	\item $W_{V_j} \subset H$, $j\in\mathbb Z_+$ is a family of sets and $F$ satisfies \VL on $W \times W_{V_j}$ for all $j\in\mathbb Z_+$,
	\item $Z\subset W$ and $T>0$ are such that for all $x_0\in Z$ the solution of the system (\ref{eq:varSys})
	 with initial conditions $x(0)=P_n x_0\in H_n$ and $V(0)=\mathrm{Id}_{H_n}$ satisfies
	\begin{eqnarray*}
	 	x^n(t)=\in W \quad \textrm{for }n>M,\ t\in[0,T],\\
		V^n_{\ast j}(t) \in W_{V_j} \quad \text{for } n>j,\ t\in[0,T].
	\end{eqnarray*}
	\end{enumerate}
	Then there exists a continuous function $\varphi:[0,T]\times Z\to W$ and
	there exists $V:[0,T] \times Z  \to \mathrm{Lin}(H,H)$, such that $V_{ij}:[0,T]\times Z \to {\mathbb R}$ for $i,j \in \mathbb{Z}$ are continuous and the following properties are satisfied.
	\begin{description}
    \item[1.]{\bf Uniform convergence:} The functions $\widehat{\varphi}^n(t,x_0) := x^n(t)$ converge uniformly to $\varphi$ on $[0,T]\times Z$.
	\item[2.] {\bf Existence and uniqueness within $W$:} For all $x\in Z$ the function $u(t):=\varphi(t,x)$ is a unique solution to the initial value problem $u'=F(u)$, $u(0)=x$ and satisfying $u(t)\in W$ for $t\in[0,T]$.
	\item[3.]  {\bf Lipschitz constant:} For $x,y\in Z$ and $t\in[0,T]$ there holds
	\begin{equation*}
	\|\varphi(t,x) - \varphi(t,y) \| \leq e^{lt}\|x - y\|.% \label{eq:LipConstant}
	\end{equation*}
	\item[4.] {\bf Convergence of variational equations:} \newline
	For each $j$ the function $\widehat{V}^n_{\ast j}(t,x_0):= \iota_n V^n_{\ast j}(t)$ converges to $V_{\ast j}$ uniformly on $[0,T] \times Z$, and
	\begin{eqnarray*}
		\|V(t,x)\| &\leq& e^{l t}, \quad (t,x) \in [0,T] \times Z,\\
		V(t,x)e_j &=& \sum_{i\in\mathbb Z_+} V_{ij}(t,x)e_i,
	\end{eqnarray*}
	and for every $a \in H$ the map $[0,T] \times Z \ni (t,x) \mapsto V(t,x)a$ is continuous.
		\item[5.] {\bf Smoothness:}  For any $x,y \in Z$ and any $t\in[0,T]$ we have
		\begin{equation*}
			\varphi(t,x) - \varphi(t,y) =  \int^1_0 V(t,y+s(x-y))ds \cdot (x - y),
			%\label{eq:intpar}
		\end{equation*}
		and  for every $j$ the partial derivative of the flow $\frac{\partial \varphi}{\partial u_j}(t,u)$ exists and
		\begin{equation*}
			\frac{\partial \varphi}{\partial u_j}(t,u)=V_{\ast j}(t,u).%  \label{eq:dfiduj=Vij}
		\end{equation*}
	\item[6.]{\bf Equation for $V$:}  $V(t,u)$ satisfies the following variational equation
	\begin{equation}
		\frac{d V_{*j}}{dt}(t,u) = \sum_i e_i \sum_k \frac{\partial F_i}{\partial u_k}(\varphi(t,u))
		V_{kj}(t,u), \label{eq:varinfdim}
	\end{equation}
	with the initial condition $V(0)=\mathrm{Id}$ in the following sense:  for each $j$ the derivative  $\frac{d V_{*j}}{dt}(t,u)$ exists, the series on r.h.s. of (\ref{eq:varinfdim}) converges uniformly
	on $[0,T] \times Z$ and equation (\ref{eq:varinfdim}) is satisfied.	
	\end{description}
\end{theorem}

Observe that the essential difficulty in application of the above theorem is to find sets $W$ and $W_{C_j}$. A systematic construction of these sets is based on the \emph{isolation property} of the vector field $F$, which will be described in Section~\ref{subsec:iso-prop}.

\subsubsection{More subtle estimates for variational matrix}

Theorem~\ref{thm:limitLN} guarantees that within sets satisfying the \emph{convergence conditions} the objects we are going to compute (or estimate) indeed exists. The rough estimates from this theorem, in particular bound on the Lipschitz constant $\|V(t,x)\|\leq e^{lt}$, are of minor practical relevance. To make our algorithm work efficiently we need  accurate bounds for block of variational matrix $V(t,x)$.

 We  consider a  decomposition of $H$ of the following form
\begin{equation}
 H=\bigoplus_{k \leq m} \langle e_k \rangle   \oplus Y, \quad Y=(I-P_m)H. \label{eq:Hdecmp}
 \end{equation}
with $m$ 1-dimensional subspaces $\langle e_i \rangle$ and infinite-dimensional subspace $Y$.

To  distinguish between $V_{ij}$, when $i,j$ are treated as indices of particular coordinates and the situation when they are supposed to denote blocks  for the block entries involving $Y$ we will write $V_{yy}$, $V_{iy}$, $V_{yi}$. In this way $V_{ij}$ for any $i,j$ will denote always a number in the matrix $V$.

 The block notation introduced above make sense also for $F^n$ (the vector field for Galerkin projection). Below we write $\frac{\partial F^n}{\partial y}$
meaning that the partial derivatives are meant with respect $y \in P_n Y$, i.e. this is a finite dimensional matrix.

The following theorem is a specialization of Theorem 15 from \cite{WZ24} to decomposition of $H$ given by (\ref{eq:Hdecmp})
\begin{theorem}{\cite[Theorem 15]{WZ24}}
\label{thm:lineqestm-inf}
Same assumptions and notation as in Theorem~\ref{thm:limitLN}. For some $m>0$ consider decomposition given by (\ref{eq:Hdecmp}).

 Assume that matrix $J \in \mathbb{R}^{(m+1) \times (m+1)}$ satisfies
\begin{equation*}%\label{eq:def-J}
 J_{ij} \geq
  \begin{cases}
     \sup_{n>m} \sup_{w \in W} \left\|\frac{\partial F^n_i}{\partial u_j}(w) \right\|,  & \text{for $i \neq j$, $i,j \leq m$}, \\
      \sup_{n>m} \sup_{w \in W}  \frac{\partial F^n_i}{\partial u_j}(w),   & \text{for $i=j$, $i\leq m$}, \\
     \sup_{n>m} \sup_{w \in W}  \mu\left(\frac{ \partial F^n_y }{\partial y}(w) \right),     & \text{for $i=j=m+1$} \\
     \sup_{n>m}   \sup_{w \in W} \left\|\frac{ \partial F^n_i}{\partial y}(w) \right\|,  & \text{for $i \leq m$, $j=m+1$} \\
       \sup_{n>m}  \sup_{w \in W} \left\| \frac{\partial  F^n_y}{\partial u_j}(w) \right\|,  & \text{for $i=m+1$, $j \leq m$}.
  \end{cases}
\end{equation*}

Then for any  $u \in Z$ and  $t \in [0,h]$
holds
\begin{eqnarray*}
  |V_{ij}(t,u)| &\leq&  \left(e^{Jt}\right)_{ij},  \quad i,j=1,\dots,m,  \\% \label{eq:Vij-estm}\\
  \|V_{yj}(t,u)\|  &\leq&  \left(e^{Jt}\right)_{m+1,j},  \quad j=1,\dots,m,\\%  \label{eq:Vyj-estm} \\
  \|V_{iy}(t,u) \| &\leq&  \left(e^{Jt}\right)_{i,m+1}, \quad i=1,\dots,m, \\% \label{eq:Vjy-estm} \\
  \|V_{yy}(t,u)\| &\leq&  \left(e^{Jt}\right)_{m+1,m+1}. % \label{eq:Vyy-estm}
\end{eqnarray*}
\end{theorem}

\subsection{Isolation property}
\label{subsec:iso-prop}

In Theorem~\ref{thm:limitLN} the crucial assumption was the existence of uniform a-priori bounds $W$ and $W_V$. It turns out that such sets can be found in algorithmic way using the \emph{isolation property} of the vector field discussed. The existence of such a-priori bounds will be later used in construction of an algorithm for rigorous integration of the flow and associated variational equations.

In what follows for  $W \subset H$ we set $W_k=\pi_k (W)$. We define the isolation property in the following way.
\begin{definition}\cite[Def. 8]{WZ24}
	Let $W\subset H$ be a set satisfying condition {\bf S} in a \gss space $H$. We say that  vector field $F:H\to H$ satisfies the \emph{isolation property} on the set $W$ if $F$ satisfies condition \textbf{C} on $W$ and there exists $K_0\in\mathbb N$ such that
	\begin{eqnarray*}
		\forall_{k\geq K_0}\, \exists_{T_k>0}\, \forall_{t\in(0,T_k]}\forall_{u\in W}
		\left(  u_k\in\mathrm{bd} W_k\Longrightarrow u_k + tF_k(u) \in\mathrm{int}W_k\right).
	\end{eqnarray*}	
\end{definition}
Geometrically the above condition means that the vector field $F$ is pointing inwards each component $W_k$ on all far coordinates of $W$. In particular, each component $W_k$, $k\geq K_0$ has nonempty interior.

\begin{definition}\cite[Def. 9]{WZ24}
	We say that a vector field $F:H\to H$ is isolating on the family of sets $\mathcal W\subset 2^H$ if $F$ satisfies the isolation property on each set $W\in\mathcal W$.
\end{definition}
Notice, that the constant $K_0$ for each set in the family $\mathcal W$ can be different.

\begin{theorem}\cite[Thm.~17 and Thm.~18]{WZ24} %\label{thm:enclosureExists}
	Assume the vector field $F:H\to H$ is isolating on the family of sets $\mathcal W$ satisfying the following properties:
	\begin{itemize}
		\item there exists $E\in \mathcal W$, such that $0\in\mathrm{int}E_k$ for all $k\in\mathbb Z_+$;
		\item $\mathcal W$ is closed with respect to the addition
		\begin{eqnarray*}
			W_1,W_2\in \mathcal W &\Longrightarrow&  W_1+W_2:=\{w_1+w_2: w_1\in W_1,\, w_2\in W_2 \}\in\mathcal W.
		\end{eqnarray*}
		\item for all $W\in\mathcal W$ the vector field $F_V$ is isolating on $W\times E$ and satisfies condition \VL
	\end{itemize}
	Then for any $Z\in\mathcal W$ there exists $T>0$ such that the assumptions of Theorem~\ref{thm:limitLN} are satisfied for the set $W=Z+E$ and $W_{V_j}=C_jE$ for some $C_j>0$, $j\in\mathbb Z_+$.
\end{theorem}

To illustrate the \emph{isolation property} let us consider a class of dissipative PDEs of the following form
\begin{equation}
	u_t = L u + N\left(u,Du,\dots,D^ru\right), \label{eq:genpde}
\end{equation}
where $u \in \mathbb{R}^n$,  $x \in \mathbb{T}^d=\left(\mathbb{R}\mod 2\pi\right)^d$, $L$ is a linear operator, $N$ is a polynomial and by $D^s u$ we denote $s^{\text{th}}$ order derivative of
$u$, i.e. the collection of all spatial partial derivatives of $u$ of
order $s$. The reason to consider polynomial and not more general functions $N$ is technical --- we need to compute
the Fourier coefficients of $N\left(u,Du,\dots,D^ru\right)$. This can be achieved by taking suitable convolutions of Fourier expansions of $u$ and its spatial partial derivatives.

We require, that the operator $L$ is diagonal in the Fourier basis
$\{e^{ikx}\}_{k \in \mathbb{Z}^d}$,
\begin{equation*}
	L e^{ikx}= -\lambda_k e^{ikx},
\end{equation*}
and there are constants  $K,L_*,L^* \geq 0$ and $p>r$ such that
\begin{equation*}
	L_* |k|^p \leq \lambda_k \leq  L^* |k|^p, \qquad \text{for all }|k| > K.
\end{equation*}

If the solutions are sufficiently smooth, the problem (\ref{eq:genpde}) can be written as an infinite ladder of ordinary differential equations for the Fourier
coefficients in $u(t,x)=\sum_{k \in \mathbb{Z}^d} u_k(t) e^{i kx}$, as follows
\begin{equation}
	\frac{d u_k}{dt} = F_k(u)=-\lambda_k u_k + N_k\left(\{u_j\}_{j \in \mathbb{Z}^d}\right), \qquad \mbox{for all
		$k \in \mathbb{Z}^d$}. \label{eq:fueq}
\end{equation}

Define a class of compact sets (in $H=c_0$ with supremum norm)
\begin{equation}\label{eq:generalFormOfW}
	W=\left\{ \{u_k\}_{k \in \mathbb{Z}^d}\,|\,  |u_k| \leq \frac{C}{|k|^s q^{|k|}}\right\}
\end{equation}
for some $q\geq 1$, $C>0$, $s>0$. Projection of the set $W$ onto $k^{\mathrm{th}}$ mode is a closed disc (or interval) centred at zero and of radius $r_k =\frac{C}{|k|^s q^{|k|}}$.

In \cite{ZKS3} it has been shown that for $q=1$
\begin{equation}\label{eq:vfPointingInvards}
	\mbox{if} \quad u \in W, \ |u_k|=\frac{C}{|k|^s}, \qquad \mbox{then} \qquad  u_k \cdot F_k(u) <0.
\end{equation}
Geometrically (\ref{eq:vfPointingInvards}) means that if $|u_k|=r_k$ for some $k> K$ then the $k$-th component of the vector field is pointing inwards the set. As a consequence, the only way a trajectory may leave the set $W$ forward in time is by increasing some of leading modes $|u_k|$, $k\leq K$ above the threshold $r_k$. However, one has to be careful with what "inwards" means in the above statement because $W$ has empty interior as a compact subset in infinite dimensional space.

From (\ref{eq:vfPointingInvards}) we can infer, that the vector field (\ref{eq:fueq}) is isolating on the family of sets (\ref{eq:generalFormOfW}) with $q=1$ -- see also \cite[Lem.~31 and Lem.~32]{WZ24}.

The isolation property is used in our approach to obtain a priori bounds for $u_k(h)$ for small $h>0$ and $|k|>K$, while the leading modes $u_k$ for $|k| \leq K$ are computed using  tools for rigorous integration of finite dimensional ODEs \cite{CAPDREVIEW,Lo,NJP} based on the interval arithmetics \cite{Mo}, i.e. to obtain sets $W$ and $W_{C_j}$ as in Theorem~\ref{thm:limitLN},  for details see section 4.4.2 in \cite{WZ}.   Moreover, from the point of view of topological method, the isolation property is of crucial importance as it shows that on the tail we have the entry behaviour, which enable us to apply the finite dimensional tools from the dynamics - the covering relations discussed in Section~\ref{sec:geometricTools}. 

\subsection{Convergence conditions and the isolation property for the KS equation}

In this work, we would like to apply Theorem~\ref{thm:limitLN} to the KS equation (\ref{eq:fuKS}). Following \cite{WZ}, we use sets of geometrically decaying tails (\ref{eq:geometricSetForm}) to integrate the KS equation and associated variational equations, which are a special case of (\ref{eq:generalFormOfW}) with $s=0$ and $q>1$.

In \cite[Section 7.6, Section 7.7]{WZ24} it is shown that the convergence conditions \textbf{C} and \VL are satisfied for a general class of vector field (\ref{eq:fueq}) on the torus and on the sets of the form (\ref{eq:generalFormOfW}) with polynomial decay, that is with $q=1$. Since each set of the form (\ref{eq:generalFormOfW}) can be included in a set with polynomial decay of coefficients, the convergence conditions \textbf{C} and \VL are also satisfied on sets (\ref{eq:generalFormOfW}). In particular, this is true for sets with geometrically decaying coefficients (\ref{eq:geometricSetForm}) and for the class of systems (\ref{eq:fueq}) on the torus.

Observe also, that if a vector field $F:H\to H$ satisfies convergence conditions \textbf{C} and \VL on some class of sets, then the same holds true for $F$ restricted to an invariant subspace of $H$.

Summarizing, the KS equation with periodic and odd boundary conditions (\ref{eq:fuKS}) satisfies the convergence conditions \textbf{C} and \VL on the sets $\WqS{q}{S}$ and $\WqS{q}{S} \times \WqS{q}{S_c}$, respectively.

The isolation property on the sets $\WqS{q}{S}$ does not follow directly from the same property for sets with polynomial decay. We will now prove that it also holds such sets

\begin{lemma}
Consider equation (\ref{eq:fuKS}) and the induced variational system. Then isolation property holds on the sets $\WqS{q}{S}$ and $\WqS{q}{S} \times \WqS{q}{S_c}$ for these equations.
\end{lemma}
\textbf{Proof:}
 To simplify  notation we set $L_k=k^2(1-\nu k^2)$. Then, for $a\in \WqS{S_E}{q}$, $S_E>0$ using (\ref{eq:fuKS}) we have
\begin{eqnarray*}
	a_kF_k(a)\leq L_k a_k^2 + k(k-1)S_E^3q^{-2k} + 2kS_E^3q^{-2k}\sum_{n=1}^\infty q^{-2n} \\
	= L_k a_k^2 + k(k-1)S_E^3q^{-2k} + 2kS_E^3q^{-2k}(q^2-1)^{-1}.
\end{eqnarray*}
Thus, if $|a_k|=S_Eq^{-k}$ then
\begin{eqnarray*}
	a_kF_k(a)\leq S_E^2q^{-2k}(L_k + O(k^2))<0
\end{eqnarray*}
for large $k$ because $L_k$ is growing to minus infinity like $-k^4$. From the above estimates we see, that for all large enough $n\geq K$  the vector field in $k$-th direction
points inside $\pi_k \WqS{q}{S}$ which proves that the isolation property is satisfied.

Now, let us consider the variational part. From (\ref{eq:fuKS}) we see that
\begin{eqnarray*} % \label{eq:partialDerKS}
	\frac{\partial F_k}{\partial a_i}(a)=   \left\{
	\begin{array}{ll}
		2k \omega(k,i) a_{|k-i|} + 2k a_{k+i}, & \hbox{ for $i \neq k$;} \\
		k^2(1-\nu k^2) +  2k a_{2k}, & \hbox{for $i=k$;} \\
	\end{array}
	\right.,
\end{eqnarray*}
where $\omega(k,i)=-1$ if $k>i$ and $\omega(k,i)=1$ if $k<i$. Thus, for $a\in \WqS{S_E}{q}$ and $C\in \WqS{S_V}{q}$ $S_E, S_V>0$ we have
\begin{equation*}
  \left| \frac{\partial F_k}{\partial a_i}(a) \right|  \leq 4k S_E q^{-|k-i|}  ,\quad \mbox{for $i \neq k$}
\end{equation*}
and we obtain a bound
\begin{eqnarray*}
	\widetilde{DF_k}(a)C&\leq& L_k C_k + 2kS_ES_Vq^{-3k} + 4kS_ES_Vq^{-k}\sum_{i=1, i\neq k}^\infty q^{-|k-i|} \\
	&\leq& L_k C_k + 2kS_ES_Vq^{-3k} + 8kS_ES_Vq^{-k}\sum_{i=1}^\infty q^{-i}\\
	&=& L_k C_k + 2kS_ES_Vq^{-3k} + 8kS_ES_Vq^{-k}(q-1)^{-1}.
\end{eqnarray*}
Now, if $|C_k|=S_Vq^{-k}$ then  we obtain
\begin{eqnarray*}
	C_k\widetilde{DF_k}(a)C&\leq& L_k S_V^2q^{-2k} + 2kS_ES_V^2q^{-4k} + 8kS_ES_V^2q^{-2k}(q-1)^{-1}\\
	&=& S_V^2q^{-2k}\left(L_k + 2kS_Eq^{-2k} + 8kS_E(q-1)^{-1}\right)\\
	&=& S_V^2q^{-2k}\left(L_k + O(k)\right).
\end{eqnarray*}
Again, for large $k$ the above expression is negative which proves the isolation property for the variational part. 
\qed

%% file: c1alg.tex
\section{The $\mathcal C^1$ algorithm for integration of dissipative PDEs.}
\label{sec:c1algorithm}

There are several algorithms in the literature \cite{C,ZKS3,ZKSper,WZ} that compute bounds of the solutions on dissipative PDEs of the form (\ref{eq:genpde}). By a $\mathcal C^1$ algorithm we understand computation of bounds on the solutions to the initial value problem  (\ref{eq:c0IVP})-(\ref{eq:c1IVP}).

Formally speaking $DF(x(t))$ in (\ref{eq:c1IVP}) may not exist on $H$, but its domain is dense in $H$. Moreover, if we restrict initial conditions to sets possessing good convergence properties (so that the assumptions of Theorem~\ref{thm:limitLN} are satisfied), then equation (\ref{eq:c1IVP}) makes sense and has a unique solution.

For finite dimensional ODEs very efficient $\mathcal C^1$ algorithms \cite{Z02,WW} are known along with their publicly available implementations \cite{CAPDREVIEW}. To the best of our knowledge, there is only one $\mathcal C^1$ algorithm \cite{AK10} that computes bounds of (\ref{eq:c1IVP}) in infinite dimensional case. It this section we propose a new algorithm which simultaneously computes bounds for the solutions (\ref{eq:c0IVP}) and associated variational equations (\ref{eq:c1IVP}).

\subsection{An outline  of $\mathcal C^1$-integration algorithm}
%\label{sec:outline-C1}
In the subsequent discussion we will quite often drop the dependence of $V$ on $u_0$, i.e. we will write just $V(t)$ instead of $V(t,u_0)$ when $u_0$ is known from the context.

Each column $V_{*j}(t)$ solves the following non-autonomous linear equation (we use decomposition $F=L+N$ as in (\ref{eq:fueq}))
\begin{equation}\label{eq:vareq}
    \frac{d}{dt}V_{*j}(t) = \frac{\partial F}{\partial u}(\varphi(t,u_0)) V_{*j}(t) = L V_{*j}(t) + \frac{\partial  N}{\partial u}(\varphi(t,u_0)) V_{*j}(t),
\end{equation}
but for different initial conditions -- that is $V_{*j}(0) = (V_0)_{*j}$.

The structure of equation (\ref{eq:vareq}) is similar to that of (\ref{eq:fueq}) -- the linear part $L$ is the same and we expect that it should dominate nonlinear part on higher modes under the mild assumptions  on the nonlinear part $N$. Thus, in principle any $\mathcal C^0$ algorithm can be applied to compute finite number of columns $V_{*j}(t)$ for $j=1,\dots,m$. The remaining block of infinite number of columns will be bounded uniformly by an operator norm.

For our further considerations it is convenient to split the phase space into $u=(x,y)$, where $x\in \mathbb{R}^m=:X$, $y$ is in the complementary subspace $Y$ and use the notation
\begin{equation*}
  V(t) =\left[\begin{array}{cc}
                  V_{xx}(t) & V_{xy}(t) \\
                   V_{yx}(t) & V_{yy}(t)
               \end{array}
        \right].
\end{equation*}
The norms in $X$ and $Y$ is quite arbitrary.  We just demand that $H=X \oplus Y$ is \gss and norms on $X$ and $Y$ are inherited from $H$.

In what follows we assume that the blocks $V_{xx}(t)$ and $V_{yx}(t)$ are computed using our $\mathcal C^0$ algorithm for integration of (\ref{eq:fueq}) (see \cite{WZ}).
Therefore for $V_{ij}(t)$ for $i,j=1,\dots,m$ are given coordinate-wise interval bounds and for  $V_{ij}$ with $j \leq m$ and $i>m$ we will have uniform bounds decaying with $i \to \infty$ depending on some constants computed by the algorithm, just like in $\mathcal C^0$ computation of $u(t)$.

Blocks  $V_{xy}$ and $V_{yy}$  will be represented in our algorithm by bounds of their operator norms. These norms will be controlled by means of the estimates from Theorem~\ref{thm:lineqestm-inf}.

Finally, we will discuss an important part of the algorithm, that is propagation of products in $V_{xx}$ and $V_{xy}$, which helps in reducing overestimation caused by the wrapping effect \cite{Lo}.

\subsection{Input and output of one step of integration scheme.}
As already mentioned, our $\mathcal C^1$ algorithm builds on the $\mathcal C^0$ algorithm proposed in \cite{WZ}. Inputs to the $\mathcal{C}^0$ algorithms which computes bounds for the main equation of (\ref{eq:c0IVP}) are

\begin{itemize}
\item $[u](t_0)$ -- a set of initial conditions,
\item $h_0>0$ -- a candidate for the time step.
\end{itemize}
Additional inputs to the $\mathcal C^1$ algorithm (integration of (\ref{eq:c1IVP})) are
\begin{itemize}
    \item $[V_{*j}](t_0)$, $j=1,\ldots,m$ -- leading columns of $[V](t_0)$, which is an initial condition for variational equation,
    \item $z(t_0)\in\mathbb R^{m+1}$ -- vector of bounds on operator norms of rows in $V_{xy}(t)$, that is $\|\left(V_{xy}\right)_{i*}(t_0)\|\leq z_i(t_0)$, for $i=1,\ldots, m$ and of
    $\|V_{yy}(t_0)\|\leq z_{m+1}(t_0)$.
\end{itemize}
Successfully completed one step of integration scheme returns the following

\noindent\textbf{Output:}
\begin{itemize}
    \item $h\leq h_0$ -- accepted time step,
    \item $[u](t_0+h)$ -- a bound for the solutions at time $t_0+h$,
    \item $[V_{*j}](t_0+h)$, $j=1,\ldots,m$ -- a bound on leading columns of $[V](t_0+h)$,
    \item $z(t_0+h)\in\mathbb R^{m+1}$ -- such that $\|\left(V_{xy}\right)_{i*}(t_0+h)\|\leq z_i(t_0+h)$, for $i=1,\ldots, m$ and
$\|V_{yy}(t_0+h)\|\leq z_{m+1}(t_0+h)$,
    \item $[E]$ -- a rough enclosure, that is a bound for $[u](t)\subset [E]$ for $t\in[t_0,t_0+h]$, such that the vector field satisfies convergence condition \textbf{C} on $[E]$,
    \item $[E_V]$ -- a rough enclosure for variational equation, that is a bound for $[V](t)\subset [E_V]$, for $t\in[t_0,t_0+h]$, such that the vector field satisfies convergence conditions \VL on $[E]\times[(E_{V})_{\ast j}]$ for $j\in\mathbb Z_+$.
\end{itemize}

Note, that usually we are interested in relatively short integration time, that is computation of $[u](T)$ for  $T>0$ being a  return to the Poincar\'e section. For this purpose the algorithm is called iteratively. The enclosures $[E]$ and $[E_V]$ are needed to compute intersection of a set of trajectories with surfaces, such as in the case of finding bounds on Poincar\'e maps and their derivatives -- see Section~\ref{sec:boundsPoincareMap}.

The very preliminary step of the $\mathcal C^0$ algorithm (see \cite{Z02,WZ}) is computation of the so called \emph{rough enclosure}. Given a set of initial conditions $[u](t_0)$ and a trial time step $h_0$ the routine returns a pair
\begin{equation}\label{eq:encC0}
\left([E],h\right)\gets\texttt{enclosure}([u](t_0),h_0),
\end{equation}
where $h\leq h_0$ is an accepted time step and $[u](t)\subset [E]$ for $t\in[t_0,t_0+h]$. In particular, this routine validates that the solutions starting from $[u](t_0)$ are defined for $t\in[t_0,t_0+h]$. If this step fails, the algorithm returns \texttt{Failure} and we stop the computation.

In what follows we will show how to compute an enclosure $[E_V]$. We will use Theorem~\ref{thm:lineqestm-inf} for the decomposition $\bigoplus_{i=1}^m \langle e_i \rangle  \oplus Y$, i.e. $X$ is decomposed into $1$-D subspaces
given by coordinate directions and initial condition $W(0)=V(t_0)$.

Since the vector field is an explicit expression which satisfies on our enclosures the convergence conditions {\rm \textbf{C}} and \VL, we can compute a matrix $J\in\mathbb R^{(m+1)\times(m+1)}$, such that for $u\in [E]$ and $A:=A(u)=DF(u)$ there holds (component-wise inequalities)
\begin{equation*}
J\geq \begin{bmatrix}
A_{11} & |A_{12}| & \cdots& |A_{1m}| & \left\|\left[A_{1y}\right]\right\|\\
|A_{21}| & A_{22} & \cdots& |A_{2m}| & \left\|\left[A_{2y}\right]\right\|\\
\vdots & \vdots & \ddots & \vdots & \vdots\\
|A_{m1}| & \cdots & |A_{m,m-1}| &A_{mm} & \left\|\left[A_{my}\right]\right\|\\
\left\|\left[A_{y1}\right]\right\| & \cdots & \left\|\left[A_{y,m-1}\right]\right\| & \left\|\left[A_{ym}\right]\right\| & \mu\left(\left[A_{yy}\right]\right)
\end{bmatrix},
\end{equation*}
where the operator norms depend on the choice of norms in $X$ and $Y$.

Since  $\left\|[V_{xy}]_{i*}(t_0)\right\|\leq z_i(t_0)$ for $i=1,\ldots,m$ and $\left\|[V_{yy}](t_0)\right\|\leq z_{m+1}(t_0)$, by Theorem~\ref{thm:lineqestm-inf} for $t\in[t_0,t_0+h]$ there holds
\begin{equation}\label{eq:estimatesC1y}
\begin{array}{rcl}
\left\|[V_{xy}]_{i*}(t)\right\|&\leq& \left(e^{J(t-t_0)}z(t_0)\right)_i,\quad i=1,\ldots,m,\\
\left\|[V_{yy}](t)\right\|&\leq& \left(e^{J(t-t_0)}z(t_0)\right)_{m+1}.\\
\end{array}
\end{equation}

From these estimates we can obtain bounds
\begin{equation}\label{eq:encC1y}
\begin{array}{rcl}
    [E_V]_{xy} &=& \left\{M\in\mathbb R^m\times \mathbb R^\infty : \|M_{i*}\|\leq \left(e^{J[0,h]}z(t_0)\right)_i,i=1,\ldots,m\right\},\\
    \mbox{}[E_V]_{yy} &=& \left\{M\in\mathbb R^\infty\times \mathbb R^\infty : \|M\|\leq \left(e^{J[0,h]}z(t_0)\right)_{m+1}\right\}.
\end{array}
\end{equation}

Using the same matrix $J$ we could obtain bounds on $W_{xx}$ and $W_{yx}$ but they are too crude. The computation of these entries is discussed in Section~\ref{subsec:re-near-col}.

\subsection{Enclosure for $V_{xx}$ and $V_{yx}$.}
\label{subsec:re-near-col}
Propagation of blocks $V_{xx}$ and $V_{yx}$ is more subtle and before we give a routine for computation of rough enclosure let us motivate why we need to compute several objects.

It is well known that integration of differential equations may suffer from the wrapping effect \cite{Lo}. In order to reduce this overestimation Lohner \cite{Lo} proposes propagation of coordinate systems along the trajectory. This means, that typical $\mathcal C^0$ solver is almost $\mathcal C^1$ as it keeps tract of approximate solution to variational equations. This observation lead to the $\mathcal C^1$ algorithm for finite dimensional systems \cite{Z02}. In our algorithm we take the advantage from the linearity of the variational equation, that is its product property
$$
V(t+h,u_0,V_0) = V(h,u_0(t),\mathrm{Id})\cdot V(t_0,u_0,V_0),
$$
where by $V(t,u_0,V_0)$ we understand a solution to variational equations after time $t$ and with the initial conditions $u_0$ and $V_0$ for the main equation and for the variational equation, respectively.

In the infinite-dimensional case, however, this approach leads to mixing of ``near'' and ``far'' columns in evaluation of this product, which causes huge overestimation. This is due to the fact, that ``far'' objects are roughly bounded by their norms, while we would like to keep explicit intervals in the finite dimensional block $V_{xx}$. In order to avoid this we propose to split the initial condition for variational equations into two parts -- as described in the sequel.

Now, let us denote by $V(t,u,V_{xx},V_{xy},V_{yx},V_{yy})$ a solution to variational equation (\ref{eq:c1IVP}) with initial condition $u$ for the $\mathcal C^0$ part and $V=\begin{bmatrix} V_{xx} & V_{xy}\\ V_{yx} & V_{yy}\end{bmatrix}$ for the $\mathcal C^1$ part. To make the notation less cumbersome we will keep only one argument in  block component of $V$, that is
$$
\begin{bmatrix}V_{xx}(t) & V_{xy}(t) \\ V_{yx}(t)&V_{yy}(t)\end{bmatrix}
:= V(t,u,V_{xx},V_{xy},V_{yx},V_{yy}).
$$

Given that the equation for variational equation is linear and non-au\-to\-nomous, we have
\begin{equation*}
V(t,u,V_{xx},V_{xy},V_{yx},V_{yy}) = V(t,u,V_{xx},0,V_{yx},0) + V(t,u,0,V_{xy},0,V_{yy})
\end{equation*}
and $V_{xx}(t), V_{yx}(t)$ appear in $V(t,u,V_{xx},0,V_{yx},0)$, only. Since
\begin{eqnarray*}
    V(t,u,V_{xx},V_{xy},V_{yx},V_{yy}) = V(t,u,\mathrm{Id},0,0,\mathrm{Id}) \begin{bmatrix}V_{xx} & V_{xy} \\ V_{yx}&V_{yy}\end{bmatrix},
\end{eqnarray*}
we have
\begin{eqnarray*}
   V(t,u,V_{xx},0,V_{yx},0) &=& V(t,u,\mathrm{Id},0,0,\mathrm{Id}) \begin{bmatrix}V_{xx} & 0 \\ 0 &0\end{bmatrix} \\
    & & +  V(t,u,\mathrm{Id},0,0,\mathrm{Id}) \begin{bmatrix}0 & 0 \\ V_{yx}& 0\end{bmatrix}.
\end{eqnarray*}
It is easy to see that
\begin{eqnarray*}
V(t,u,\mathrm{Id},0,0,\mathrm{Id}) \begin{bmatrix}V_{xx} & 0 \\ 0 &0\end{bmatrix}= V(t,u,\mathrm{Id},0,0,0) \begin{bmatrix}V_{xx} & 0 \\ 0 &0\end{bmatrix}
\end{eqnarray*}
and
\begin{eqnarray*}
V(t,u,\mathrm{Id},0,0,\mathrm{Id}) \begin{bmatrix}0 & 0 \\ V_{yx}& 0\end{bmatrix} =  V(t,u,0,0,V_{yx},0)  %\label{eq:var-ev-Vyx}.
\end{eqnarray*}
This gives
\begin{eqnarray}
\begin{bmatrix}V_{xx}(t) & 0 \\ V_{yx}(t)&0\end{bmatrix} = V(t,u,\mathrm{Id},0,0,0)\begin{bmatrix}V_{xx} & 0 \\ 0&0\end{bmatrix}  + V(t,u,0,0,V_{yx},0). \label{eq:progVxxVyx}
\end{eqnarray}

Observe that in the product of matrices in (\ref{eq:progVxxVyx}) we do not mix ``near'' and ``far'' columns of $V$ because the second matrix contains only one nonzero block $V_{xx}$. The second component $V(t,u,0,0,V_{yx},0)$ is computed independently and we expect it should be very small due to the strong dissipation on higher modes.

Summarizing, in order to efficiently propagate bounds for $m$ leading columns of $V$ we propose  in each time step to integrate two different initial conditions. Denote them by
\begin{equation*}
\widetilde{V} = \begin{bmatrix}\mathrm{Id} & 0 \\ 0 & 0\end{bmatrix},\qquad
\widehat{V} = \begin{bmatrix}0 & 0 \\ [V_{yx}](t_0) & 0\end{bmatrix}.
\end{equation*}

For each of the two initial conditions the rough enclosure for $m$ leading columns of the solution is computed by the same routine and for the same vector field but for different initial conditions. Note, that in order to bound vector field for variational equation we need already computed bound $[E]$ for the main equation (see (\ref{eq:encC0})). We will emphasize this adding extra argument $[E]$ in the following sequence of calls:
\begin{equation*}
\begin{array}{rcl}
\left([\widetilde E_V]_{*1},h\right)&\gets&\texttt{enclosure}([E],\widetilde V_{*1},h),\\
&\ldots&\\
\left([\widetilde E_V]_{*m},h\right)&\gets&\texttt{enclosure}([E],\widetilde V_{*m},h),\\
\left([\widehat E_V]_{*1},h\right)&\gets&\texttt{enclosure}([E],\widehat V_{*1},h),\\
&\ldots&\\
\left([\widehat E_V]_{*m},h\right)&\gets&\texttt{enclosure}([E],\widehat V_{*m},h).
\end{array}
\end{equation*}
Note, that each call to \texttt{enclosure} can adjust the time step. If none of the above calls returns \texttt{Failure}, then we have computed a rough enclosure over time step $h$ for the main equation and for $m$ leading columns of the variational equation for two different initial conditions.

Now, let us set
\begin{equation}\label{eq:encC1x}
\begin{array}{rcl}
    [E_V]_{xx} &=& [\widetilde E_V]_{xx} [V_{xx}](t_0) + [\widehat E_V]_{xx},\\
    \mbox{}[E_V]_{yx} &=& [\widetilde E_V]_{yx} [V_{xx}](t_0) + [\widehat E_V]_{yx}.
\end{array}
\end{equation}
Gathering (\ref{eq:encC0}), (\ref{eq:encC1y}) and (\ref{eq:encC1x}) we obtain the following partial output of our algorithm.

\textbf{Partial output:}
\begin{itemize}
    \item $h$ -- adjusted and finally accepted time step,
    \item $[E]$ -- rough enclosure for the main equation,
    \item $[E_V]$ -- rough enclosure for the variational equation.
\end{itemize}

\subsection{Propagation of $V(t)$.}

Following \cite{Z02} we represent finite block of $V$ as the doubleton
$$[V_{xx}](t_0) = D + C[R] + B[R_0],$$
where $D,C,B$ are matrices and $[R], [R_0]$ are interval matrices containing zero.

The propagation of this block is realized as (compare (\ref{eq:progVxxVyx}))
\begin{eqnarray}
 [V_{xx}](t_0+h) &:=& \pi_{xx}V\left(h,[u](t_0),\mathrm{Id},0,0,0\right) (D + C[R] + B[R_0]) \nonumber \\
  &+ & \pi_{xx}V(h,[u](t_0),0,0,[V_{yx}](t_0),0).  \label{eq:newVxx}
\end{eqnarray}
In order to reduce the wrapping effect the product of matrices in (\ref{eq:newVxx}) is propagated as suggested by Lohner \cite{Lo} (see also \cite{Z02,WW}). After that, a rather thin part (i.e. having a small diameter) $\pi_{xx}V(h,[u](t_0),0,0,[V_{yx}](t_0),0)$ is absorbed into $B[R_0]$ in representation of $[V_{xx}](t_0+h)$ -- see \cite[Evaluation 4 in Section 3.2]{Z02}.

Each column of $V_{yx}$ is represented as a geometric series with decay $q>1$ and a constant $C_j=C_j(t_0)$, that is
$$
|V_{ij}(t_0)|\leq C_jq^{-i},\quad i>m,\ j=1,\ldots,m.
$$
We would like to compute new constants $C_j(t_0+h)$, $j=1,\ldots,m$, so that
\begin{equation}\label{eq:newVyx}
|V_{ij}(t_0+h)|\leq C_j(t_0+h)q^{-i},\quad i>m,\ j=1,\ldots,m.
\end{equation}

From (\ref{eq:progVxxVyx}) we obtain
\begin{eqnarray*}
[V_{yx}](t_0+h) &=& \pi_{yx}V\left(h,[u](t_0),\mathrm{Id},0,0,0\right) [V_{xx}](t_0) \nonumber \\
&+ & \pi_{yx}V(h,[u](t_0),0,0,[V_{yx}](t_0),0),
\end{eqnarray*}
which for $j=1,\ldots,m$ gives
\begin{equation}
C_j(t_0+h) = \sum_{k=1}^m \widetilde{C_k} \left|[V_{kj}](t_0)\right| + \widehat{C_j},  \label{eq:newCjVyx}
\end{equation}
where $\widetilde{C_j}$ and $\widehat{C_j}$ are computable constants satisfying
\begin{eqnarray*}
|V_{ij}(h,[u](t_0),\mathrm{Id},0,0,0)|&\leq& \widetilde{C_j}q^{-i},\\
|V_{ij}(h,[u](t_0),0,0,[V_{yx}](t_0),0)|&\leq& \widehat{C_j}q^{-i},
\end{eqnarray*}
for $i>m,\ j=1,\ldots,m$. Indeed, these constants are computed by means of the $\mathcal C^0$ algorithm from \cite{WZ} applied to variational equation with initial conditions $\widetilde V_{*j}$ and $\widehat V_{*j}$, $j=1,\ldots,m$.

Finally, using estimates from (\ref{eq:estimatesC1y}) we can set $z(t_0+h):=e^{Jh}z(t_0)$ obtaining bounds
\begin{equation}\label{eq:newVyxVyy}
\begin{array}{rcl}
z_i(t_0+h) &\geq& \left\|[V_{xy}]_{i*}(t_0+h)\right\|,\qquad i=1,\ldots,m,\\
z_{m+1}(t_0+h)&\geq&\left\|[V_{yy}](t_0+h)\right\|.
\end{array}
\end{equation}
Gathering (\ref{eq:newVxx}), (\ref{eq:newVyx}), (\ref{eq:newCjVyx}) and (\ref{eq:newVyxVyy}) we obtain remaining output of the algorithm.

\noindent\textbf{Partial output:}
\begin{itemize}
    \item $z(t_0+h)$ -- a vector of bounds on norms in blocks $V_{xy}$ and $V_{yy}$ at time $t_0+h$,
    \item $[V_{xx}](t_0+h)$ -- finite block in the doubleton representation,
    \item $[V_{yx}](t_0+h)$ -- tails of leading $m$ columns represented as geometric series.
\end{itemize}

\subsection{Bounds on derivatives of Poincar\'e map}\label{sec:boundsPoincareMap}
Efficient computation of Poincar\'e map is performed as has been discussed in \cite{WZ,KaWZ}. In this section we show how to obtain bounds on derivatives of Poincar\'e map from derivatives of the flow.

Let $\alpha\colon H\to \mathbb R$ be a smooth map such that
\begin{equation}\label{eq:nablaalpha}
\alpha(x)=\tilde{\alpha}(x_1,x_2,\ldots,x_m),
\end{equation}
so that is partial derivatives with respect to all variables $x_j$, $j>m$ are zero. Let us define a Poincar\'e section by
$$
\Pi = \left\{x : \alpha(x) = 0\right\}.
$$
Denote by $T:H \supset \mbox{dom}(T) \to\mathbb R$  the return time to the section and by $\PM$ the Poincar\'e map for a semiflow $\varphi$ induced by $x'=f(x)$. In what follows we will derive estimates for derivatives of $T$ and $\PM$.

Differentiating the identity $\alpha(\varphi(T(x),x))=0$ with respect to $x_j$  we have
\begin{equation*}%\label{eq:DtEquation}
D\alpha(\PM(x))\left( f(\PM(x))\right) \cdot \frac{\partial T}{\partial {x_j}}(x) + D\alpha(\PM(x))\left(\frac{\partial \varphi}{\partial {x_j}}(T(x),x)\right) = 0.
\end{equation*}
The quantity
$$g:= \left(D\alpha(\PM(x))\left(f(\PM(x))\right)\right)^{-1} = \left(\sum_{k=1}^m\alpha_{x_k}(\PM(x))f_k(\PM(x))\right)^{-1}$$
is independent from $j$ and easily computable (finite sum) provided the vector field is transverse to the section at $\PM(x)$. Thus, due to (\ref{eq:nablaalpha}), partial derivatives of $T$ are given by
\begin{equation}\label{eq:dt}
 \frac{\partial T}{\partial {x_j}}(x) = -g\sum_{k=1}^m \frac{\partial \alpha}{\partial {x_k}}(\PM(x))\frac{\partial \varphi_k}{\partial x_j}(T(x),x).
\end{equation}

Differentiation of $\PM(x) = \varphi(T(x),x)$ yields
\begin{equation}\label{eq:DP}
\frac{\partial \PM_k}{\partial x_j}(x) = \frac{\partial \varphi_k}{\partial x_j}(T(x),x) + f_k(\PM(x)) \frac{\partial T}{\partial {x_j}}(x).
\end{equation}

Let us recall, that the already proposed $\mathcal C^1$ algorithm computes $\frac{\partial \varphi_k}{\partial x_j}$ as explicit intervals, for $k,j\leq m$, while three remaining blocks of $D_x\varphi$ are estimated uniformly by an operator norm -- the details on evaluation of (\ref{eq:dt}) and (\ref{eq:DP}) will be given in Appendix~\ref{app:partialDerviatives}.

%% file: appendix-num.tex
\section{Explicit formulas}

In this entire section the phase space $H=c_0$ - the sequences converging to zero,  is decomposed as $\mathcal H = X\oplus Y$ equipped with the norm $\|(x,y)\| = \max\left\{\|x\|_2,\|y\|_\infty\right\}$. For $z=(x,y)$ we also use indexing $x=(x_1,\ldots,x_m)\in X,y=(y_{m+1},\ldots)\in Y$.

Moreover, we assume that the $z,\PM(z)\in W_{q,S}$, for some constants $q>1$ and $S\geq 0 $.

\subsection{Some formulas for norms of blocks of derivatives}

%\label{app:norms-blocks-der}

 $A \in \mbox{Lin}(H,H)$. $A$ is split into four blocks denoted by
$$
A=\begin{bmatrix}
A_{xx} & A_{xy} \\ A_{yx} & A_{yy}
\end{bmatrix}.
$$
We will also use $A_{ky}$ to denote an operator $A_{ky}\colon Y\to \mathbb R$ given as the row matrix $A_{ky}=\begin{bmatrix}A_{k,m+1} & A_{k,m+2} &\cdots\end{bmatrix}$. The norms of the blocks of $A$ are given by the following explicit formulas.

\begin{eqnarray}
\|A_{ky}\|  &=&  \sum_{j>m} \left|A_{kj}\right|,  \label{eq:appDfkDy}\\
\|A_{yy}\|  &=& \sup_{k>m} \|A_{ky}\| = \sup_{k>m} \sum_{j>m} \left|A_{kj}\right|, \nonumber\\% \label{eq:appDfydy}\\
\|A_{yk}\| &=& \sup_{j >m} \left|A_{jk} \right|, \nonumber\\ %\label{eq:normDyDxk}\\
\|A_{kx}\|  &=& \left( \sum_{j\leq m} \left(A_{kj}\right)^2 \right)^{1/2} \leq  \sum_{j\leq m} \left|A_{kj}\right|,\nonumber\\
\|A_{yx}\|&=& \sup_{k>m} \left( \sum_{j\leq m} \left|A_{kj} \right|^2  \right)^{1/2} \leq
            \sup_{k>m}  \sum_{j\leq m} \left|A_{kj} \right|. \label{eq:normPyx}
\end{eqnarray}

\subsection{Partial derivatives of Poincar\'e map}
\label{app:partialDerviatives}

In this section we describe how to evaluate (\ref{eq:dt}--\ref{eq:DP}). Let us denote
\begin{equation}
W_{kj}(x)=f_k(\PM(x))\frac{\partial T}{\partial x_j}(x)  \label{eq:W-corr}
\end{equation}
In the sequel we will use the following notation $V(x)=\frac{\partial \varphi}{\partial x}(T(x),x)$ and $V_{**}$ are blocks in $V$. Then
\begin{equation*}
D\PM=V + W.
\end{equation*}

Let $D$ be a constant such that for $k>m$ there holds
\begin{equation}
|f_k(\PM(z)|\leq D\label{eq:fk-estm}
\end{equation}

Let us consider four cases corresponding to four block of $D\PM$ we have to compute.

\noindent
\textbf{Block $W_{xx}$.} For $k,j\leq m$ all terms in (\ref{eq:dt}) are given by explicit numbers (intervals).

\noindent
\textbf{Block $W_{yx}$.} Put $W_{yx}:=W_{k,j}$ for ${j\leq m<k}$. From (\ref{eq:normPyx}) and
(\ref{eq:W-corr}) we obtain
\begin{eqnarray}
    \|W_{yx}\| &\leq& \sup_{k>m}\left(\sum_{j=1}^m\left|W_{kj}\right|^2\right)^{1/2} \leq  \sup_{k>m}\sum_{j=1}^m\left|W_{kj}\right| \nonumber \\
    &\leq&  \sup_{k>m} \sum_{j=1}^m\left|f_k(\PM(z))\frac{\partial T}{\partial x_j}(z)\right|
    \leq  \max_{k>m} \left|f_k(\PM(z))\right|\sum_{j=1}^m\left|\frac{\partial T}{\partial x_j}(z)\right|, \label{eq:Pyx1}
\end{eqnarray}

 For $j\leq m$ derivatives of $\frac{\partial T}{\partial x_j}$ in (\ref{eq:dt}) are given by explicit numbers. From (\ref{eq:fk-estm})
 we obtain
\begin{equation}\label{eq:Pyx2}
\max_{k>m} \left|f_k(\PM(z))\right|\sum_{j=1}^m\left|\frac{\partial T}{\partial x_j}(z)\right|
\leq   D\sum_{j=1}^m\left|\frac{\partial T}{\partial x_j}(z)\right|.
\end{equation}
Gathering (\ref{eq:Pyx1})--(\ref{eq:Pyx2}) we get a final bound
$$
 \|W_{yx}\| \leq  D\sum_{j=1}^m\left|\frac{\partial T}{\partial x_j}(z)\right|.
$$

\noindent
\textbf{Rows $W_{ky}$.} In this case we want to estimate operator norm of each row (that is fixed $k\leq m$) separately.
 Computation similar to (\ref{eq:Pyx1}) leads to an estimate
\begin{multline}\label{eq:Pxy1}
\|W_{ky}\| \leq \left(\sum_{j>m}\left|W_{kj}\right|^2\right)^{1/2} \leq \sum_{j>m}\left|W_{kj}\right|
=  \sum_{j>m}\left|f_k(\PM(z))\frac{\partial T}{\partial x_j}(z)\right|\\
\leq  \left|f_k(\PM(z))\right|\sum_{j>m}\left|\frac{\partial T}{\partial x_j}(z)\right|.
\end{multline}
From (\ref{eq:dt}) we have (see (\ref{eq:appDfkDy}))
\begin{multline}\label{eq:Pxy2}
\sum_{j>m}\left|\frac{\partial T}{\partial x_j}(z)\right| \leq |g|\sum_{j>m}\sum_{i=1}^m\left|\alpha_{x_i}(\PM(z))\frac{\partial \varphi_i}{\partial x_j}(T(z),z)\right|\\
= |g|\sum_{i=1}^m\left(\left|\alpha_{x_i}(\PM(z))\right|\sum_{j>m}\left|\frac{\partial \varphi_i}{\partial x_j}(T(z),z)\right|\right)
= |g| \sum_{i=1}^m\left|\alpha_{x_i}(\PM(z))\right|\|V_{iy}\|.
\end{multline}
Gathering (\ref{eq:Pxy1})--(\ref{eq:Pxy2}) we get
$$
\|W_{ky}\|\leq |g|\left|f_k(\PM(z))\right|\sum_{i=1}^m\left|\alpha_{x_i}(\PM(z))\right|\|V_{iy}\|.
$$

\textbf{Block $W_{yy}$.} Computation similar to (\ref{eq:Pxy1})--(\ref{eq:Pxy2}) gives an estimate
\begin{multline*}%\label{eq:Pyy1}
\|W_{yy}\|\leq\max_{k>m}\left(|g|\left|f_k(\PM(z))\right|\sum_{i=1}^m\left|\alpha_{x_i}(\PM(z))\right|\|V_{iy}\|\right) \\ \leq |g|D\sum_{i=1}^m\left|\alpha_{x_i}(\PM(z))\right|\|V_{iy}\|.
\end{multline*}

\section{How to verify the cone condition for h-sets?}
%\label{subsubsec:cones-ver}

Recall the phase space is decomposed as $\mathcal H = X\times Y$ and equipped with the norm $\|(x,y)\| = \max\left\{\|x\|_2,\|y\|_\infty\right\}$. Let $N,M$ be $h$-sets with cones of the form
\begin{equation}
Q_{K}(x,y)= \sum_{ 1 \leq k \leq m} Q_{K,1}(x) - \|y\|_\infty^2, \label{eq:Q-diag}
\end{equation}
where
\begin{equation*}
Q_{K,1}(x)= \sum_{1 \leq k \leq m} q_{K,k} x_k^2.
\end{equation*}
and $K=N$ or $K=M$. We will write $q_k$ instead if $K$ will be clear from the context.

Let $f$ be a continuous map expressed in the coordinate systems of $N$ and $M$, respectively. Assume that the map $f(x,y)=(f_x(x,y),f_y(x,y))$ has all continuous partial derivatives and denote
$$
Df(z) = \begin{bmatrix}
\frac{\partial f_i}{\partial z_j}(z)
\end{bmatrix}_{i,j\in\mathbb N} = \begin{bmatrix} f_{xx}(z) & f_{xy}(z) \\ f_{yx}(z) & f_{yy}(z)\end{bmatrix}
$$
for $z\in N$. The existence of continuous partial derivatives allows us to use Euler formula
\begin{equation*}%\label{eq:EulerFormula}
f(z_1)-f(z_2) = \left(\int_0^1 Df(z_1+t(z_2-z_1))dt\right) (z_1-z_2).
\end{equation*}

We would like to derive set of computable inequalities that will guarantee that the cone condition is satisfied for the covering relation $N\cover{f}M$. To simplify notation will use $\|\cdot\|$ to denote norms of $x$ and $y$.

Let us fix  $z_i=(x_i,y_i)\in N$, $i=1,2$ and set $dx=x_1-x_2$ and $dy=y_1-y_2$.
We define $\overline{f_{**}}$ as the integral average of partial derivatives on interval joining $z_1$ and $z_2$, for example
\begin{equation*}
\overline{f_{xx}} = \int_0^1 f_{xx}(z_1 + t(z_2-z_1)dt.
\end{equation*}
In a similar way we define $\overline{f_{xk}}$ and $\overline{f_{ky}}$ for $k=1,\ldots,m$. Using this notation we have
\begin{equation}\label{eq:average}
\begin{array}{rcl}
f_x(z_1)-f_x(z_2) &=& \overline{f_{xx}}dx + \overline{f_{xy}}dy,\\
f_y(z_1)-f_y(z_2) &=& \overline{f_{yx}}dx + \overline{f_{yy}}dy.
\end{array}
\end{equation}
Now we derive a formula for
\begin{eqnarray*}
    Q_M(f(z_1)-f(z_2))= Q_{M,1}(f_x(z_1)-f_x(z_2)) - \|f_y(z_1)-f_y(z_2)\|^2.
\end{eqnarray*}
In what follows we Set $q_k=q_{M,k}$. Using (\ref{eq:average}) we have
\begin{eqnarray*}
    Q_{M,1}(f_x(z_1)-f_x(z_2))=Q_{M,1}(\overline{f_{xx}} dx + \overline{f_{xy}}dy), \\
    = \sum_{k\leq m} q_k \left(\overline{f_{kx}} dx + \overline{f_{ky}} dy\right)^2 = \sum_{k \leq m} q_k \left(\overline{f_{kx}} dx\right)^2 \\
    + \sum_{k \leq m} 2 q_k \left(\overline{f_{kx}} dx\right) \left(\overline{f_{ky}} dy\right) + \sum_{k\leq m} q_k \left(\overline{f_{ky}} dy\right)^2.
\end{eqnarray*}
The three expressions in the above sum can be bounded as follows
\begin{eqnarray*}
    \sum_{k \leq m} q_k \left(\overline{f_{kx}} dx\right)^2 &=& dx^T \left((\overline{f_{xx}})^TQ_{M,1} \overline{f_{xx}} \right) dx, \\
    \left|\sum_{k \leq m} 2 q_k \left(\overline{f_{kx}} dx\right) \left(\overline{f_{ky}} dy\right) \right| &\leq&  2\left(\sum_{k \leq m} |q_k| \, \|\overline{f_{kx}}\| \, \|\overline{f_{ky}}\| \right) \|dx\| \, \|dy\|, \\
    \left|\sum_{k\leq m} q_k \left(\overline{f_{ky}} dy\right)^2  \right| &\leq& \left(\sum_{k\leq m} |q_k| \cdot \|\overline{f_{ky}}\|^2\right) \|dy\|^2.
\end{eqnarray*}
Since
\begin{eqnarray*}
    \|f_y(z_1)-f_y(z_2)\| =\|\overline{f_{yx}} dx + \overline{f_{yy}} dy\| \leq  \|\overline{f_{yx}}\| \, \|dx\| + \|\overline{f_{yy}}\| \, \|dy\|
\end{eqnarray*}
we finally obtain
\begin{eqnarray*}
    Q_M(f(z_1)-f(z_2)) \geq dx^T \left((\overline{f_{xx}})^TQ_{M,1} \overline{f_{xx}} \right) dx + \\
    - 2\left(\sum_{k \leq m} |q_k| \cdot \|\overline{f_{kx}}\| \cdot \|\overline{f_{ky}}\| \right) \|dx\| \cdot \|dy\| + \\
    - \left(\sum_{k\leq m} |q_k| \cdot \|\overline{f_{ky}}\|^2\right) \|dy\|^2
    -  (\|\overline{f_{yx}}\| \cdot \|dx\| + \|\overline{f_{yy}}\| \cdot \|dy\|)^2 \\
    = dx^T \left((\overline{f_{xx}})^TQ_{M,1} \overline{f_{xx}} \right) dx - \|\overline{f_{yx}}\|^2 \|dx\|^2  + \\
    - 2\left(\left(\sum_{k \leq m} |q_k| \cdot \|\overline{f_{kx}}\| \cdot \|\overline{f_{ky}}\| \right)  + \|\overline{f_{yx}}\|  \cdot \|\overline{f_{yy}}\|  \right) \|dx\| \cdot \|dy\| + \\
    - \left( \left(\sum_{k\leq m} |q_k| \cdot \|\overline{f_{ky}}\|^2\right)  +  \|\overline{f_{yy}}\|^2 \right) \|dy\|^2.
\end{eqnarray*}

Summarizing, we obtain
\begin{equation*}
Q_M(f(z_1)-f(z_2)) - Q_N(z_1-z_2) \geq a \|dx\|^2 - 2c \|dx\| \cdot \|dy\| + d \|dy\|^2
\end{equation*}
where $a,c,d$ are such that
\begin{eqnarray*}
    dx^T \left((\overline{f_{xx}})^TQ_{M,1} \overline{f_{xx}} \right) dx - \|\overline{f_{yx}}\|^2 \|dx\|^2 - Q_{N,1}(dx) \geq a \|dx\|^2, \\
    \left(\left(\sum_{k \leq m} |q_k| \cdot \|\overline{f_{kx}}\| \cdot \|\overline{f_{ky}}\| \right)  + \|\overline{f_{yx}}\|  \cdot \|\overline{f_{yy}}\|  \right) \leq c, \\
    1 - \left( \left(\sum_{k\leq m} |q_k| \cdot \|\overline{f_{ky}}\|^2\right)  +  \|\overline{f_{yy}}\|^2 \right) > d.
\end{eqnarray*}
Observe that $a \|dx\|^2 - 2c \|dx\| \cdot \|dy\| + d \|dy\|^2$ is positive for all $dx\neq 0$ or $dy\neq 0$ if and only if $a>0$ and $ad>c^2$. Therefore we have proved the following lemma.

\begin{lemma} \label{lem:rigConeVerify}
    Assume $(N,Q_N)$, $(M,Q_M)$ are $h$-sets with cones of the form (\ref{eq:Q-diag}). Denote by $\left[Df_{xx}(N)\right]_I$ a convex hull of the set of matrices $Df_{xx}(z)$, $z\in N$.

    Assume that there is $a>0$ such that for $A\in\left[Df_{xx}(N)\right]_I$
         \begin{eqnarray*}
        \left(A^TQ_{M,1} A - Q_{N,1} - \sup_{z\in N}\|f_{yx}(z)\|^2 \mathrm{Id}\right)  - a \mathrm{Id}
    \end{eqnarray*}
    is positive definite.

    Assume also that $c\geq 0,d>0$ are constants such that
     \begin{eqnarray*}
\sup_{z\in N}\left(\left(\sum_{k \leq m} |q_{M,k}| \cdot \|f_{kx}(z)\| \cdot \|f_{ky}(z)\| \right)  + \|f_{yx}(z)\|  \cdot \|f_{yy}(z)\|  \right) \leq c, \\
         1 - \sup_{z\in N}\left( \left(\sum_{k\leq m} |q_{M,k}| \cdot \|f_{ky}(z)\|^2\right)  +  \|f_{yy}(z)\|^2 \right) > d.
     \end{eqnarray*}

    If $ ad > c^2$, then  $Q_M(f(z_2)-f(z_1))>Q_N(z_2-z_1)$ for $z_1,z_2\in N$, $z_1\neq z_2$.
\end{lemma}

\section{Changing coordinate system on the section.}
Here we want to represent $DP$ in a linear coordinate system $\widetilde A(x,y) = (A(x),y)$, that is we want to compute
$$
    M = \left(\widetilde A\right)^{-1}DP \widetilde A = \begin{bmatrix} M_{xx} & M_{xy} \\ M_{yx} & M_{yy}\end{bmatrix}.
$$
Straightforward computation gives
\begin{eqnarray*}
    M_{xx} &=& A^{-1}P_{xx}A\\
    M_{xy} &=& A^{-1}P_{xy}\\
    M_{yx} &=& P_{yx}A\\
    M_{yy} &=& P_{yy}
\end{eqnarray*}
Note, that coefficients in $P_{xx}$ are given as explicit intervals, all columns in the block $P_{yx}$ are given as geometric series and we have only bounds on the operator norms $D_yP_k$, $k=1,\ldots,n$ and $P_{yy}$. Clearly $M_{xx}$ can be computed by direct multiplication of interval matrices and $\|M_{yy}\|=\|P_{yy}\|$. For the remaining blocks we have the following estimates.

 \subsection{Computation of $M_{xy}$.}
 Put $B:=A^{-1}$ and denote by $M_{k}$ the operator norm of $k$-th row of $M_{xy}$. The range is one-dimensional and in the domain we have $l_\infty$ norm.
\begin{eqnarray*}
 M_{k} &=& \sum_{j>m}\left|\sum_{i=1}^m B_{ki}P_{ij}\right|\leq \sum_{j> m}\sum_{i=1}^m |B_{ki}||P_{ij}|= \sum_{i=1}^m |B_{ki}|\sum_{j> m}|P_{ij}|\\
  &=&
 \sum_{i=1}^m |B_{ki}|\|D_y P_i\|.
 \end{eqnarray*}

 \subsection{Computation of $M_{yx}$.}
 Here we use direct bound
 \begin{equation*}
 \|M_{yx}\|\leq \|P_{yx}\|\cdot \|A\|_2
 \end{equation*}